\documentclass[final]{siamltex}

\usepackage{amsmath,amssymb}
\usepackage{mathrsfs}
\usepackage{mathtools}
\usepackage{bm}
\usepackage{graphicx}
\usepackage{booktabs}
\usepackage{cite}
\usepackage{hyperref}
\usepackage{cleveref}
\usepackage{aliascnt}
\usepackage{enumitem}
\usepackage{xcolor}
\usepackage{multirow}
\usepackage{placeins}

\hypersetup{hidelinks}


\DeclareMathOperator{\Var}{Var}

\newaliascnt{remark}{theorem}
\newtheorem{remark}[remark]{Remark}
\aliascntresetthe{remark}
\crefname{remark}{Remark}{Remarks}
\Crefname{remark}{Remark}{Remarks}

\newaliascnt{example}{theorem}

\aliascntresetthe{example}
\crefname{example}{Example}{Examples}
\Crefname{example}{Example}{Examples}

\newaliascnt{assumption}{theorem}
\newtheorem{assumption}[assumption]{Assumption}
\aliascntresetthe{assumption}
\crefname{assumption}{Assumption}{Assumptions}
\Crefname{assumption}{Assumption}{Assumptions}

\makeatletter
\def\@begintheorem#1#2{\par\medskip\bgroup{\scshape #1\ #2. }\it\ignorespaces}
\def\@opargbegintheorem#1#2#3{\par\medskip\bgroup%
    {\scshape #1\ #2\ ({\upshape #3}). }\it\ignorespaces}
\def\@endtheorem{\egroup\par\medskip}
\makeatother

\newcommand{\R}{\mathbb{R}}
\newcommand{\E}{\mathbb{E}}
\newcommand{\N}{\mathcal{N}}
\newcommand{\bx}{\bm{x}}
\newcommand{\bz}{\bm{z}}
\newcommand{\bxi}{\bm{\xi}}

\newcommand{\email}[1]{\texttt{#1}}

\begin{document}

\title{Variance stable Gaussian smoothing estimator\\ for arbitrary high-order derivatives}

\author{Shuixin Fang\textsuperscript{\dag}
    \and Xuerui Wang\textsuperscript{\dag}
    \and Hang Zeng\textsuperscript{\dag}
    \and Tao Zhou\textsuperscript{\dag,\,*}}

\maketitle
\begingroup
    \renewcommand{\thefootnote}{\fnsymbol{footnote}}%
    \footnotetext[2]{SKLMS \& Institute of Computational Mathematics and Scientific/Engineering Computing, Academy of Mathematics and Systems Science, Chinese Academy of Sciences, Beijing 100190, China.}%
    \footnotetext[1]{Corresponding author. (\email{tzhou@lsec.cc.ac.cn}).}%
\endgroup


\begin{abstract}
    Computing high-order partial derivatives in high dimensions by automatic differentiation can be costly in both computation and memory.
    Stein's identity offers an alternative through Gaussian smoothing, but the variance of the resulting estimators can diverge as the step size tends to zero, even with antithetic sampling.
    We propose a variance stabilization technique that extends antithetic sampling from a reflection pair to $n$ equally spaced points on a circle in the complex plane.
    Under local analyticity, we prove that, when $n$ is at least the derivative order, the bias is of order $\mathcal{O}(\sigma^n)$ and the variance remains bounded as the smoothing step size $\sigma$ tends to zero.
    Numerical experiments support these predictions and demonstrate the effectiveness of the estimator in physics-informed neural networks (PINNs) for polyharmonic and Cahn--Hilliard equations.
    The technique mitigates the bias-variance trade-off in smoothing-based derivative estimation and is useful for PINNs solving high-order partial differential equations.
\end{abstract}



\pagestyle{myheadings}
\thispagestyle{plain}


\section{Introduction}\label{sec_intro}

Physics-informed neural networks (PINNs) have become a popular deep learning framework in scientific computing, owing to their mesh-free formulation, robust approximation, rapid convergence, and strong generalization capabilities~\cite{Raissi2019,Karniadakis2021,Jagtap2020Adaptive,Kawaguchi2016,Kawaguchi2022Generalization}.
Despite their success, PINNs can incur substantial computational costs when applied to PDEs with high input dimensions, high derivative orders, or both.
A major bottleneck is the evaluation of network derivatives by automatic differentiation, whose computational and memory costs can grow rapidly as either the input dimension or the derivative order increases~\cite{Hu2024STDE}.
Both types of difficulty arise in scientific and engineering applications.
For example, high-dimensional PDEs arise as Hamilton--Jacobi--Bellman equations in stochastic optimal control and nonlinear Black--Scholes equations in multi-asset derivative pricing~\cite{Han2018DeepBSDE}, while high-order PDEs model thin-film dynamics and crystal growth~\cite{Oron1997ThinFilm,Elder2002PFC}.
Efficient derivative computation is therefore essential for extending the applicability of PINNs to high-dimensional and high-order PDEs, motivating a series of recent works~\cite{He2023,Hu2024SDGD,Hu2024HTE,Lyu2022MIM,Hu2025Bias,Cai2025,Cai2026DeepMartNet,Han2020DerivativeFree}.

Randomized-smoothing PINNs (RS-PINNs) offer a promising approach to reducing derivative evaluation costs~\cite{He2023,Hu2025Bias}.
By perturbing network inputs with Gaussian noise, they represent the approximate solution as a Gaussian-smoothed neural network.
Stein's identity then expresses the input derivatives as expectations involving network evaluations rather than explicit derivatives~\cite{Fournie1999Malliavin,Stein1981}.
These expectations can be approximated via Monte Carlo sampling, bypassing the nested automatic differentiation required to evaluate PDE residuals.
Repeated application of Stein's identity extends this formulation to higher-order derivatives, making it an attractive approach for high-dimensional and high-order PDEs.

Despite these computational advantages, Monte Carlo approximation introduces new challenges in training RS-PINNs.
First, nonlinear combinations of Monte Carlo estimates in the loss function introduce gradient bias during training with stochastic gradient descent.
Hu et al.~\cite{Hu2025Bias} analyzed this bias and developed correction techniques, showing that debiasing can increase sampling variance and computational cost.
Second, high-order derivative estimation poses a distinct challenge: variance divergence as the smoothing step size vanishes.
While a small smoothing step size is desirable for preserving the network's approximation capability, it causes the variance of standard Stein estimators to diverge.
Variance stabilization has been achieved for low-order derivatives.
For first-order derivatives, Jongeneel et al.~\cite{JongeneelYueKuhn2024} developed a complex-step randomized estimator whose variance remains bounded as the step size vanishes.
Antithetic sampling stabilizes second-order derivative estimators~\cite{BalasubramanianGhadimi2022,He2023}, but it does not generally prevent variance divergence for higher-order derivatives.
This limitation motivates the development of variance-stable estimators for arbitrarily high-order derivatives.

To overcome this limitation, we extend antithetic sampling~\cite{HammersleyMorton1956} from a symmetric pair on the real line to rotationally symmetric points in the complex plane.
Specifically, for each Gaussian direction, we average the Stein estimators over $n$ equally spaced complex increments on a circle, yielding the so-called Rotationally Symmetric Stein Estimator (RSSE).
When $n$ is at least the target derivative order, this rotational averaging eliminates the lower-order contributions responsible for variance divergence.
We establish two theoretical properties of RSSE under suitable analyticity and integrability assumptions and demonstrate its numerical advantages in PINNs:
\begin{itemize}
    \item \textbf{Variance stability.} We rigorously prove its variance converges to a finite limit as the step size $\sigma$ vanishes.
    \item \textbf{High-order bias decay.} We also prove that its approximation bias is of order $\mathcal{O}(\sigma^n)$.
    \item \textbf{Efficient PINN solver.} Combining RSSE with PINNs yields efficient solvers for high-order and high-dimensional PDEs.
    Numerical experiments demonstrate that the method outperforms antithetic Stein estimators across most test cases.
\end{itemize}

The remainder of this paper is organized as follows.
Section~\ref{sec_prelim} reviews the notation, Gaussian smoothing identities, and multivariate Hermite polynomials.
Section~\ref{sec_formulation} develops the RSSE and describes its use in PINNs.
Section~\ref{sec_thm} establishes theoretical bounds on the approximation bias and variance.
Section~\ref{sec_experiments} reports numerical results for derivative approximation and PDE solution.
Section~\ref{sec_concl} concludes the paper.

\section{Preliminaries}\label{sec_prelim}

This section introduces the notation, reviews Gaussian smoothing and Stein's identity, and recalls multivariate Hermite polynomials and their orthogonality relations.

\subsection{Notation}\label{sec_notation}

We introduce some notation used throughout this paper.
Let $\mathbb{Z} := \{\ldots, -2, -1, 0, 1, 2, \ldots\}$, $\mathbb{N}_0 := \{0, 1, 2, \ldots\}$, $\mathbb{N}_ + := \{1, 2, 3, \ldots\}$, and $\mathbb{C}^d$ be the $d$-dimensional complex space.
For any set $E$ with finitely many elements, denote by $|E|$ the number of elements in $E$.
Define $\mathbf{e}_i = (0, \ldots, 0, 1, 0, \ldots, 0)^{\top} \in \R^d$ with the $1$ in the $i$-th position.
Denote by $\mathrm{Tr}[M]$ the trace of a square matrix $M$.
Denote by $\mathcal{M}$ the set of multi-indices $\alpha = (\alpha_1, \ldots, \alpha_d) \in \mathbb{N}_0^d$ of length $d$.
For every $\alpha \in \mathcal{M}$ and every $\bm{z} \in \mathbb{C}^d$, we write
\begin{equation*}
    |\alpha| := \sum_{i = 1}^{d} \alpha_i, \quad \alpha! := \prod_{i = 1}^{d} \alpha_i!, \quad \bm{z}^{\alpha} := \prod_{i = 1}^{d} z_i^{\alpha_i}, \quad \partial^{\alpha} := \partial_{x_1}^{\alpha_1} \cdots \partial_{x_d}^{\alpha_d}.
\end{equation*}
For $i \in \mathbb{N}_0$, the sum $\sum_{|\alpha| = i}$ is taken over all $\alpha \in \mathcal{M}$ with $|\alpha| = i$.
We write $\alpha = \beta$ for $\alpha, \beta \in \mathcal{M}$ if and only if $\alpha_i = \beta_i$ for all $i = 1, \ldots, d$.

Throughout this paper, $\bm{\xi}$ denotes a standard Gaussian random vector, that is, $\bm{\xi} \sim \mathcal{N}(0, I_d)$.
For every complex-valued random variable $Z$, denote its variance by $\Var[Z] := \E[|Z - \E[Z]|^2]$.
For a sequence of random variables $(X_{\sigma})_{\sigma > 0}$ and a sequence of positive real numbers $(z_{\sigma})_{\sigma > 0}$, write $X_{\sigma} = \mathcal{O}_p(z_{\sigma})$ as $\sigma \to 0$ if for every $\varepsilon > 0$, there exist $C_{\varepsilon}, \sigma_{\varepsilon} > 0$, both independent of $\sigma$, such that $\mathbb{P}(|X_{\sigma}| > C_{\varepsilon} z_{\sigma}) < \varepsilon$ for every $\sigma \in (0, \sigma_{\varepsilon})$.
For a random event $A$, we denote its indicator function by $\mathbf{1}_A$, namely, $\mathbf{1}_A = 1$ if $A$ occurs and $\mathbf{1}_A = 0$ otherwise.
For a complex number $z \in \mathbb{C}$, we denote its complex conjugate by $\overline{z}$, and denote its real and imaginary parts by $\mathrm{Re}(z)$ and $\mathrm{Im}(z)$, respectively.

\subsection{Gaussian smoothing and derivative estimation}\label{sec_rs}

For a real-valued function $u: \R^d \to \R$ and a step size $\sigma > 0$, define its Gaussian smoothing by
\begin{equation}\label{eq_smoothed}
    u_{\sigma}(\bx) := \E\bigl[u(\bx + \sigma\bxi)\bigr], \quad \bxi \sim \mathcal{N}(0, I_d).
\end{equation}
Under sufficient regularity and integrability to justify differentiation and Gaussian integration by parts, Stein's identity~\cite{Stein1981} gives
\begin{equation}\label{eq_stein_lap}
    \Delta u_{\sigma}(\bx) = \E\!\left[\frac{|\bxi|^2 - d}{\sigma^2}\,u(\bx + \sigma\bxi)\right].
\end{equation}
A sample average of the integrand is unbiased for $\Delta u_{\sigma}(\bx)$.
As an approximation to $\Delta u(\bx)$, it also incurs a smoothing bias, generally of order $\mathcal{O}(\sigma^2)$ under appropriate smoothness and moment conditions.
Repeated integration by parts gives analogous representations for higher-order and mixed derivatives.

The factors $\sigma^{-p}$ appearing in estimators of derivatives of order $p$ can amplify lower-order Taylor terms and produce large sampling variance at small step sizes.
For the Laplacian, symmetry of the Gaussian law and $\E[|\bxi|^2 - d] = 0$ yield the centered antithetic representation
\begin{equation}\label{eq_anti_lap}
    \Delta u_{\sigma}(\bx) = \E\!\left[\frac{|\bxi|^2 - d}{2\sigma^2} \bigl(u(\bx + \sigma\bxi) + u(\bx - \sigma\bxi) - 2u(\bx)\bigr)\right].
\end{equation}
The centered increment is of order $\sigma^2$, so its division by $\sigma^2$ removes the small-step singularity when the corresponding second moments are controlled.
For higher-order derivatives, a reflection pair does not generally cancel all lower-order terms.
The rotational averaging introduced in Section~\ref{sec_formulation} provides the additional cancellations needed for arbitrary target orders.

\subsection{Hermite polynomials}\label{sec_hermite}

We now briefly review the multivariate Hermite polynomials, a fundamental tool underlying the RSSE.
For every multi-index $\alpha \in \mathcal{M}$, the multivariate Hermite polynomial $H_{\alpha}: \mathbb{C}^d \to \mathbb{C}$ is defined by
\begin{equation}\label{eq_hermite}
    H_{\alpha}(\bx) := (-1)^{|\alpha|} \rho(\bx)^{-1} \partial^{\alpha} \rho(\bx), \quad \rho(\bx) := (2\pi)^{-d / 2} \exp(-\bx^{\top}\bx / 2), \quad \bx \in \mathbb{C}^d.
\end{equation}
The definition in \eqref{eq_hermite} is the tensor-product, unnormalized version of the probabilists' Hermite Rodrigues formula; see \cite[Definition~3 and (6)]{Rahman2017WienerHermite}.
Repeated Gaussian integration by parts, equivalently the tensor-product orthogonality in \cite[Corollary~10]{Rahman2017WienerHermite}, then gives
\begin{equation}\label{eq_ortho}
    \E[H_{\alpha}(\bxi) \bxi^{\beta}] = \alpha! \delta_{\alpha\beta}, \quad \forall \alpha, \beta \in \mathcal{M} \text{ with } |\beta| \le |\alpha|,
\end{equation}
where $\delta_{\alpha\beta} := 1$ if $\alpha = \beta$, and $\delta_{\alpha\beta} := 0$ otherwise.

\section{Method}\label{sec_formulation}

In this section, we present the core formulation of the proposed RSSE.
We first derive a single-directional Monte Carlo estimator for arbitrary high-order derivatives, and then stabilize its variance by averaging over rotationally symmetric increments in the complex plane, yielding the RSSE.
Finally, we describe its application to PINNs.

To present the method, we consider a function $u: \mathbb{R}^d \to \mathbb{R}$ and aim to approximate its derivative $\partial^{\alpha} u(\bx)$ at a given point $\bx \in \R^d$ for a multi-index $\alpha \in \mathcal{M}$ with $|\alpha| \ge 1$.
To avoid technicalities, we initially assume that $u$ admits an analytic extension to $\mathbb{C}^d$ and that all infinite sums and expectations in this section can be interchanged.
In fact, our method only requires local analyticity around $\bx$ together with a bounded measurable extension; this setting is treated rigorously in Section~\ref{sec_thm}.

\subsection{Basic Expansion}\label{sec_basic}

We first extend the real step size $\sigma$ in \eqref{eq_stein_lap} to a complex step size $z \in \mathbb{C}$.
Specifically, the Taylor expansion of $u(\bx + z \bxi)$ around $z = 0$ yields
\begin{equation}\label{eq_dir_taylor}
    u(\bx + z \bxi) - u(\bx) = \sum_{m = 1}^{\infty} z^m \sum_{|\beta| = m} \frac{\bxi^{\beta}}{\beta!}\,\partial^{\beta} u(\bx), \quad z \in \mathbb{C}.
\end{equation}
Multiplying both sides of \eqref{eq_dir_taylor} by the Hermite polynomial $H_{\alpha}(\bxi)$ introduced in \eqref{eq_hermite} and taking expectations with respect to $\bxi$ gives
\begin{align}
    \E\Bigl[H_{\alpha}(\bxi) \bigl(u(\bx + z \bxi) - u(\bx) \bigr)\Bigr] = \;& \sum_{m = 0}^{|\alpha|} z^m \sum_{|\beta| = m} \frac{\partial^{\beta} u(\bx)}{\beta!} \alpha! \delta_{\alpha\beta} + \sum_{m = |\alpha| + 1}^{\infty} z^m R_m(\alpha) \notag\\
    = \;& z^{|\alpha|} \partial^{\alpha} u(\bx) + \sum_{m = |\alpha| + 1}^{\infty} z^m R_m(\alpha)\label{eq_extract}
\end{align}
for every $z \in \mathbb{C}$, where the first equality follows from \eqref{eq_ortho}, and $R_m(\alpha)$ is the remainder term independent of $z$, defined by
\begin{equation}\label{eq_remd}
    R_m(\alpha) := \sum_{|\beta| = m} \frac{\partial^{\beta} u(\bx)}{\beta!} \E[H_{\alpha}(\bxi) \bxi^{\beta}].
\end{equation}

For $z \ne 0$, rearranging \eqref{eq_extract} yields a probabilistic representation of the target derivative $\partial^{\alpha} u(\bx)$:
\begin{equation}\label{eq_prim_mc}
    \partial^{\alpha} u(\bx) - \E\bigl[\mathcal{D}_{z, \bxi}^{\alpha} u(\bx)\bigr] = -\sum_{m = |\alpha| + 1}^{\infty} z^{m - |\alpha|} R_m(\alpha), \quad z \in \mathbb{C} \setminus \{0\},
\end{equation}
where $\mathcal{D}_{z, \bxi}^{\alpha} u(\bx)$ is the single-directional estimator defined by
\begin{equation}\label{eq_sing_stein}
    \mathcal{D}_{z, \bxi}^{\alpha} u(\bx) := z^{-|\alpha|} H_{\alpha}(\bxi) \bigl(u(\bx + z \bxi) - u(\bx)\bigr).
\end{equation}
The expansion~\eqref{eq_prim_mc} provides an approximation of $\partial^{\alpha} u(\bx)$ by approximating the expectation with Monte Carlo samples of $\bxi$, where Gaussian parity eliminates $R_{|\alpha| + 1}(\alpha)$, so the bias is $\mathcal{O}(|z|^2)$ when the displayed expansion is justified.
However, the variance of this estimator can diverge as $z \to 0$ for $|\alpha| \ge 2$; see Remark~\ref{rmk_var_div} below.
Section~\ref{sec_rsse} addresses this issue.

\begin{remark}[Divergence of variance for single-directional estimator]\label{rmk_var_div}
    For simplicity, we assume the remainder on the right-hand side of \eqref{eq_prim_mc} is of order $\mathcal{O}(|z|)$ as $z \to 0$, which is commonly satisfied.
    Then the variance of $\mathcal{D}_{z, \bxi}^{\alpha} u(\bx)$ admits the following expansion as $z \to 0$:
    \begin{align}
        \Var\bigl[\mathcal{D}_{z, \bxi}^{\alpha} u(\bx)\bigr] & = \E\bigl[\bigl|\mathcal{D}_{z, \bxi}^{\alpha} u(\bx)\bigr|^2\bigr] - \bigl|\E\bigl[\mathcal{D}_{z, \bxi}^{\alpha} u(\bx)\bigr]\bigr|^2 \notag\\
        & = |z|^{-2 |\alpha|} \E\bigl[H_{\alpha}(\bxi)^2 \bigl|u(\bx + z \bxi) - u(\bx)\bigr|^2\bigr] - \bigl|\partial^{\alpha} u(\bx)\bigr|^2 + \mathcal{O}(|z|) \notag\\
        & = |z|^{2 - 2|\alpha|} \E\bigl[H_{\alpha}(\bxi)^2 \bigl(\bxi^{\top} \nabla u(\bx)\bigr)^2\bigr] - \bigl|\partial^{\alpha} u(\bx)\bigr|^2 + \mathcal{O}\bigl(|z|^{3 - 2|\alpha|}\bigr), \label{eq_expand_var}
    \end{align}
    where the last equality follows from Taylor's expansion of $z \mapsto u(\bx + z \bxi)$ at $z = 0$.
    Consequently, $\Var\bigl[\mathcal{D}_{z, \bxi}^{\alpha} u(\bx)\bigr]$ diverges whenever $\nabla u(\bx) \neq 0$ and $|\alpha| \ge 2$.
\end{remark}

\subsection{RSSE}\label{sec_rsse}

To obtain a variance-stable estimator, we average \eqref{eq_prim_mc} for $z$ at multiple points $\sigma_1, \sigma_2, \ldots, \sigma_n \in \mathbb{C}$, which yields
\begin{align}
    \partial^{\alpha} u(\bx) - \E\Bigl[\frac{1}{n} \sum_{j = 1}^{n} \mathcal{D}_{\sigma_j, \bxi}^{\alpha} u(\bx)\Bigr] & = -\sum_{m = |\alpha| + 1}^{\infty} \frac{1}{n} \sum_{j = 1}^{n} \sigma_j^{m - |\alpha|} R_m(\alpha). \label{eq_avg}
\end{align}
In particular, let $\sigma_1, \ldots, \sigma_n$ be equally spaced points on a circle of radius $\sigma$ in $\mathbb{C}$, i.e.,
\begin{equation}\label{eq_defsgmj}
    \sigma_j := \sigma \omega^j, \quad \omega := \exp(2\pi \mathrm{i} / n), \quad j = 1, 2, \ldots, n, \quad \mathrm{i} := \sqrt{-1}.
\end{equation}
The above choice of $\sigma_j$ leads to 
\begin{equation}\label{eq_filter}
    \frac{1}{n} \sum_{j = 1}^{n} \sigma_j^{m - |\alpha|} = \begin{cases}
        \sigma^{m - |\alpha|}, & \text{if } m \in \{|\alpha| + kn\}_{k \in \mathbb{Z}}, \\
        0, & \text{otherwise.}
    \end{cases}
\end{equation}
The expansion \eqref{eq_avg} then reduces to
\begin{equation}\label{eq_avg_rot}
    \partial^{\alpha} u(\bx) - \E\Bigl[\mathcal{E}^{\alpha}[u](\bx; \bxi)\Bigr] = -\sum_{k = 1}^{\infty} \sigma^{k n} R_{|\alpha| + k n}(\alpha) = \mathcal{O}(\sigma^n), \quad \text{as } \sigma \to 0,
\end{equation}
where $\mathcal{E}^{\alpha}[u](\bx; \bxi)$ denotes the proposed RSSE given by
\begin{equation}\label{eq_rsse}
    \mathcal{E}^{\alpha}[u](\bx; \bxi) := \frac{1}{n} \sum_{j = 1}^{n} \mathcal{D}_{\sigma_j, \bxi}^{\alpha} u(\bx) = \frac{1}{n} \sum_{j = 1}^{n} \sigma_j^{-|\alpha|} H_{\alpha}(\bxi) \bigl(u(\bx + \sigma_j \bxi) - u(\bx)\bigr)
\end{equation}
for $\bx \in \R^d$.
Here, the second equality follows from the definition of $\mathcal{D}_{\sigma_j, \bxi}^{\alpha} u(\bx)$ in \eqref{eq_sing_stein}, and the dependence of $\mathcal{E}^{\alpha}$ on $\sigma$ and $n$ is suppressed for notational brevity.

\begin{remark}[Connection to the antithetic estimator]\label{rem_antithetic}
    With $n = 2$ in \eqref{eq_rsse}, the antithetic Laplacian estimator \eqref{eq_anti_lap} proposed in \cite[(11)]{He2023} coincides exactly with $\sum_{i = 1}^d \mathcal{E}^{2\mathbf{e}_i}[u](\bx; \bxi)$, where the step size $\sigma$ is identified with the real smoothing bandwidth appearing in \eqref{eq_anti_lap}.
    Hence, the RSSE can be viewed as a natural generalization of \eqref{eq_anti_lap} to arbitrary derivative orders.
\end{remark}

\begin{remark}[Connection to classical numerical differentiation]
    The roots-of-unity filtering in~\eqref{eq_filter} is closely related to classical numerical differentiation of analytic functions via equispaced samples on a circle in the complex plane~\cite{LynessMoler1967,Bornemann2011,AustinKravanjaTrefethen2014}; see~\cite[Secs.~2.3 and~2.7]{AustinKravanjaTrefethen2014} for an analogous formulation of the corresponding Taylor-coefficient and derivative formulas.
\end{remark}

We next examine the variance for $n \ge |\alpha|$.
Since the variance is bounded above by the second moment, it suffices to establish the following bound:
\begin{align}
    \E\Bigl[\Big| \mathcal{E}^{\alpha}[u](\bx; \bxi) \Big|^2\Bigr] = \;& \E\Bigl[\Big| \frac{1}{n} \sum_{j = 1}^{n} \sigma_j^{-|\alpha|} H_{\alpha}(\bxi) \bigl(u(\bx + \sigma_j \bxi) - u(\bx)\bigr) \Big|^2\Bigr] \notag\\
    = \;& \E\Bigl[\Big| \frac{H_{\alpha}(\bxi)}{n} \sum_{j = 1}^{n} \sum_{m = 1}^{\infty} \sigma_j^{m - |\alpha|} \sum_{|\beta| = m} \frac{\bxi^{\beta}}{\beta!}\,\partial^{\beta} u(\bx) \Big|^2\Bigr] \notag\\
    = \;& \E\Bigl[\Bigl| H_{\alpha}(\bxi) \sum_{|\beta| = |\alpha|} \frac{\bxi^{\beta}}{\beta!}\,\partial^{\beta} u(\bx) \Bigr|^2\Bigr] + \mathcal{O}(\sigma^n), \label{eq_2nd_mom}
\end{align}
as $\sigma \to 0$, where the first line follows from the definition in \eqref{eq_sing_stein} and Taylor's expansion of $\sigma_j \mapsto u(\bx + \sigma_j \bxi)$ at $\sigma_j = 0$, and the last line follows from \eqref{eq_filter}.

In short, expansions \eqref{eq_avg_rot} and \eqref{eq_2nd_mom} show that the bias of $\mathcal{E}^{\alpha}[u](\bx; \bxi)$ vanishes to order $\sigma^n$, while its second moment remains bounded as $\sigma \to 0$, meaning that its variance is also bounded, in contrast to the single-directional estimator \eqref{eq_sing_stein}.
And the antithetic estimator \eqref{eq_anti_lap} enjoys this boundedness only for derivatives no higher than second order.
These observations are sharpened in Section~\ref{sec_thm} into a rigorous theorem.

Based on \eqref{eq_avg_rot}, we approximate the expectation by Monte Carlo sampling, yielding the empirical RSSE for $\partial^{\alpha} u(\bx)$:
\begin{equation}\label{eq_ksamp_rsse}
    \widehat{\partial}^{\alpha}u(\bx) := \frac{1}{K}\sum_{k = 1}^{K}\mathcal{E}^{\alpha}[u](\bx; \bxi_k) = \frac{1}{n K}\sum_{k = 1}^{K} \sum_{j = 1}^{n}\sigma_j^{-|\alpha|}H_{\alpha}(\bxi_k)\bigl(u(\bx + \sigma_j\bxi_k) - u(\bx)\bigr),
\end{equation}
where $\bxi_1, \ldots, \bxi_K$ are independent samples from $\mathcal{N}(0, I_d)$; the second equality follows from the definition in \eqref{eq_rsse}; $H_{\alpha}$ and $\sigma_j$ are defined in \eqref{eq_hermite} and \eqref{eq_defsgmj}, respectively.
The dependence on $\sigma$, $n$, and $K$ is suppressed in the notation $\widehat{\partial}^{\alpha}$.
With $n \geq |\alpha|$, we can expect the estimator \eqref{eq_ksamp_rsse} satisfies the error estimate
\begin{equation}\label{eq_mc_error}
    \widehat{\partial}^{\alpha}u(\bx) - \partial^{\alpha} u(\bx) = \mathcal{O}(\sigma^n) + \mathcal{O}_p(K^{-1 / 2})
\end{equation}
as $\sigma \to 0$ and $K \to \infty$, where the first term on the right-hand side reflects the approximation bias from \eqref{eq_avg_rot}, and the second term captures the Monte Carlo sampling error governed by the second-moment bound \eqref{eq_2nd_mom}.
Rigorous bias and variance bounds for \eqref{eq_rsse} are established in Theorem~\ref{thm_main} in the next section.

\subsubsection{Applications to PINNs}\label{sec_monte_carlo}

We now apply the RSSE to PINNs~\cite{Raissi2019} for solving the following PDE problem:
\begin{equation}\label{eq_pde}
    \begin{cases}
        \mathcal{A}u_* = f_*, & \bx \in \Omega, \\
        \mathcal{B}u_* = g_*, & \bx \in \Gamma,
    \end{cases}
\end{equation}
where $u_* = u_*(\bx) \in \R$ is the unknown solution; $f_* = f(\bx, u_*(\bx))$ and $g_* := g(\bx, u_*(\bx))$ are forcing terms with given $f, g: \R^d \times \R \to \R$; $\Omega \subset \R^d$ is a bounded domain with boundary $\Gamma := \partial\Omega$.
The differential operators $\mathcal{A}$ and $\mathcal{B}$ take the form
\begin{align}
    \mathcal{A}u(\bx) & := \sum_{|\alpha| \le p} A_{\alpha}(\bx, u(\bx))\,\partial^{\alpha}(F_{\alpha} \circ u)(\bx), \label{eq_op}\\
    \mathcal{B}u(\bx) & := \sum_{|\alpha| \le q} B_{\alpha}(\bx, u(\bx))\,\partial^{\alpha}(G_{\alpha} \circ u)(\bx), \label{eq_op_b}
\end{align}
where $p, q \in \mathbb{N}_0$ are the respective derivative orders; $A_{\alpha}, B_{\alpha}: \R^d \times \R \to \R$ and $F_{\alpha}, G_{\alpha}: \R \to \R$ are all given functions.

Applying \eqref{eq_mc_error} to the positive-order derivatives in \eqref{eq_op}--\eqref{eq_op_b} gives
\begin{align}
    \widehat{\mathcal{A}}u(\bx) & := \sum_{|\alpha| \le p} A_{\alpha}(\bx, u(\bx))\,\widehat{\partial}^{\alpha}(F_{\alpha} \circ u)(\bx), \label{eq_hatA} \\
    \widehat{\mathcal{B}}u(\bx) & := \sum_{|\alpha| \le q} B_{\alpha}(\bx, u(\bx))\,\widehat{\partial}^{\alpha}(G_{\alpha} \circ u)(\bx), \label{eq_hatB}
\end{align}
with $n \ge \max\{p, q\}$ and the convention $\widehat{\partial}^{\alpha}v := v$ in case $|\alpha| = 0$.
In the spirit of PINNs, combined with the above RSSE approximations, we approximate $u_*$ by a neural network $u_{\theta}$ with trainable parameters $\theta$, trained by minimizing the loss function
\begin{equation}\label{eq_pinn_loss}
    \begin{aligned}
        \bar{\mathcal{L}}(\theta) := \;& \frac{1}{|S(\Omega)|}\sum_{\bx \in S(\Omega)} \bigl\{(\widehat{\mathcal{A}}u_{\theta} - f_{\theta})(\bx)\bigr\} \bigl\{(\widehat{\mathcal{A}}^{\prime}u_{\theta} - f_{\theta})(\bx)\bigr\} \\
        & + \frac{\lambda}{|S(\Gamma)|}\sum_{\bx \in S(\Gamma)} \bigl\{(\widehat{\mathcal{B}}u_{\theta} - g_{\theta})(\bx)\bigr\} \bigl\{(\widehat{\mathcal{B}}^{\prime}u_{\theta} - g_{\theta})(\bx)\bigr\}.
    \end{aligned}
\end{equation}
Here, $f_{\theta}(\bx) := f(\bx, u_{\theta}(\bx))$ and $g_{\theta}(\bx) := g(\bx, u_{\theta}(\bx))$.
The sets $S(\Omega)$ and $S(\Gamma)$ consist of collocation points randomly sampled from $\Omega$ and $\Gamma$, respectively, at each training iteration, and $\lambda > 0$ is a hyperparameter balancing the PDE and boundary losses.
The operators $\widehat{\mathcal{A}}^{\prime}$ and $\widehat{\mathcal{B}}^{\prime}$ are independent copies of $\widehat{\mathcal{A}}$ and $\widehat{\mathcal{B}}$ in \eqref{eq_hatA}--\eqref{eq_hatB}, obtained by replacing $\bxi_1, \ldots, \bxi_K$ in \eqref{eq_ksamp_rsse} with their independent copies $\bxi_1^{\prime}, \ldots, \bxi_K^{\prime}$.
Using such independent copies to mitigate the bias of the loss function~\eqref{eq_pinn_loss} is a standard technique in the literature~\cite{Hu2025Bias,Guo2022fPINNs,Cai2025,cai2024martnet}.


\section{Theoretical properties}\label{sec_thm}

This section presents the theoretical analysis of the bias and variance of the RSSE~\eqref{eq_rsse} under the following assumption on $u$.

\begin{assumption}
    \label{assu_main}
    Let $\partial^{\alpha} u(\bx)$ be the derivative to be approximated, with given $\bx \in \R^d$ and $\alpha \in \mathcal{M}$ such that $|\alpha| \ge 1$.
    The function $u : \mathbb{C}^d \to \mathbb{C}$ is measurable, real-valued on $\R^d$, with $M := \sup_{\bm{z} \in \mathbb{C}^d} |u(\bm{z})| < \infty$.
    There exists $r > 0$ such that $u$ is analytic on $D_r := \{\bx + \bm{z} : \bm{z} \in \mathbb{C}^d, \, |\bm{z}| < r\}$.
\end{assumption}

Theorem~\ref{thm_main} below is the main result of this work.
It rigorously quantifies the bias and variance of the RSSE~\eqref{eq_rsse}, with explicit constants for the remainder terms.
It also shows that the bias of the RSSE is of order $\sigma^n$, whereas its leading variance term is independent of $\sigma$, thereby establishing its variance stability.

\begin{theorem}\label{thm_main_alt}
    \label{thm_main}
    Under Assumption~\ref{assu_main}, the estimator $\mathcal{E}^{\alpha}[u](\bx; \bxi)$ defined in \eqref{eq_rsse} satisfies
    \begin{align}
        \bigl|\E\bigl[\mathcal{E}^{\alpha}[u](\bx; \bxi)\bigr] - \partial^{\alpha} u(\bx)\bigr| & \le C_{\mathrm{E}}\, \sigma^n, \label{eq_thm_mean}\\
        \Bigl|\Var\bigl[\mathcal{E}^{\alpha}[u](\bx; \bxi)\bigr] - \Var\bigl[H_{\alpha}(\bxi) \sum_{|\beta| = p} \frac{\bxi^{\beta}}{\beta!}\,\partial^{\beta} u(\bx) \bigr]\Bigr| & \le C_{\mathrm{V}}\, \sigma^n \label{eq_thm_var}
    \end{align}
    for all $n \ge p := |\alpha|$, and all $\sigma \in (0, \bar{\sigma})$ with $\bar{\sigma} := r\,\min\bigl\{1 / (8\sqrt{d}), \, 1 / \bigl(4\sqrt{2p + n}\,\bigr)\bigr\}$, where $C_{\mathrm{E}}$ and $C_{\mathrm{V}}$ are the explicit constants given by \eqref{eq_cmean} and \eqref{eq_cvar}, respectively.
\end{theorem}


\begin{proof}
    For notational convenience, we define the integer $n^{\prime}$, the event $A_{\sigma}$ and the Taylor polynomial $T_m$ by
    \begin{equation}\label{eq_taylor}
        n^{\prime} := p + n, \quad A_{\sigma} := \{|\bxi| \le r / (2\sigma)\}, \quad T_m(\bx, \bz) := \sum_{|\beta| = m} \frac{\bz^{\beta}}{\beta!}\,\partial^{\beta} u(\bx),
    \end{equation}
    for $\sigma \in (0, \bar{\sigma})$, $\bz \in \R^d$, and $m \in \mathbb{N}_0$, where $p$ and $\bar{\sigma}$ are introduced in Theorem~\ref{thm_main}, and $r$ in Assumption~\ref{assu_main}; $\bxi$ refers to the Gaussian random vector appearing in $\mathcal{E}^{\alpha}[u](\bx; \bxi)$ from \eqref{eq_rsse}.
    This proof is divided into three steps.

    \subsubsection*{(i) Taylor expansion and probability bounds}

    On $A_{\sigma}$, for every $\theta \in [0, 1]$ and $j \in \{1, \ldots, n\}$, the definition of $D_r$ in Assumption~\ref{assu_main} implies that $\bx + \theta \sigma \omega^j \bxi \in \overline{D_{r / 2}}$, where $\overline{D_{r / 2}}$ denotes the closure of $D_{r / 2}$.
    Hence, Taylor's theorem applied to $s \mapsto u(\bx + s\omega^j\bxi)$ at $s = 0$ yields
    \begin{equation}\label{eq_ftaylor}
        \begin{aligned}
            &u(\bx + \sigma\omega^j\bxi) - u(\bx) \\
            = \;& \sum_{m = 1}^{n^{\prime} - 1} \sigma^m \omega^{jm}\,T_m(\bx, \bxi) + \sigma^{n^{\prime}} \omega^{jn^{\prime}} \int_0^1 \frac{(1 - \theta)^{n^{\prime} - 1}}{(n^{\prime} - 1)!}\,(\bxi \cdot \nabla)^{n^{\prime}} u(\bx + \theta \sigma \omega^j \bxi) \, d\theta,
        \end{aligned}
    \end{equation}
    for every $\sigma > 0$, on the event $A_{\sigma}$, where $T_m$ is defined in \eqref{eq_taylor}.

    Substituting the definition \eqref{eq_sing_stein} of $\mathcal{D}_{\sigma_j, \bxi}^{\alpha} u(\bx)$ into the right-hand side of \eqref{eq_rsse}, and then the expansion \eqref{eq_ftaylor} for the resulting increment $u(\bx + \sigma\omega^j\bxi) - u(\bx)$, we obtain
    \begin{equation}\label{eq0_app_dec}
        \begin{aligned}
            \mathcal{E}^{\alpha}[u] = \;& H_{\alpha}(\bxi) \sum_{m = 1}^{n^{\prime} - 1} \sigma^{m - p}\,T_m(\bx, \bxi)\,\Bigl(\frac{1}{n} \sum_{j = 1}^{n} \omega^{j(m - p)} \Bigr) \\
            & + \sigma^n H_{\alpha}(\bxi)\,\mathcal{R}(\bx, \bxi, \sigma) \,\Bigl(\frac{1}{n} \sum_{j = 1}^{n} \omega^{j(n^{\prime} - p)} \Bigr),
        \end{aligned}
    \end{equation}
    on the event $A_{\sigma}$, where $\mathcal{R}(\bx, \bxi, \sigma)$ is the remainder term defined by
    \begin{equation}\label{eq_rem_int}
        \mathcal{R}(\bx, \bxi, \sigma) := \frac{1}{n} \sum_{j = 1}^{n} \int_0^1 \frac{(1 - \theta)^{n^{\prime} - 1}}{(n^{\prime} - 1)!}\,(\bxi \cdot \nabla)^{n^{\prime}} u(\bx + \theta \sigma \omega^j \bxi) \, d\theta.
    \end{equation}
    Further using \eqref{eq_filter} to simplify the right-hand side of \eqref{eq0_app_dec}, and multiplying both sides by $\mathbf{1}_{A_{\sigma}}$ to extend the identity to the whole probability space, we obtain
    \begin{equation}\label{eq_app_dec}
        \mathcal{E}^{\alpha}[u]\,\mathbf{1}_{A_{\sigma}} = H_{\alpha}(\bxi)\,T_p(\bx, \bxi)\,\mathbf{1}_{A_{\sigma}} + \sigma^n H_{\alpha}(\bxi)\,\mathcal{R}(\bx, \bxi, \sigma)\,\mathbf{1}_{A_{\sigma}}
    \end{equation}
    for every $\sigma > 0$.
    We note that $\mathcal{R}(\bx, \bxi, \sigma)$ is defined by \eqref{eq_rem_int} only on $A_{\sigma}$; for rigor, $\mathcal{R}(\bx, \bxi, \sigma)\,\mathbf{1}_{A_{\sigma}}$ should be understood to be zero on $A_{\sigma}^c$.

    Now we bound the probability of the complement $A_{\sigma}^c$ of the good event $A_{\sigma}$, which will be used in the remaining steps to bound the terms in \eqref{eq_app_dec}.
    By \cite[Lemma~1]{LaurentMassart2000}, the distribution of $\bxi \sim \mathcal{N}(0, I_d)$ implies
    \begin{equation}\label{eq0_lm}
        \mathbb{P}\bigl(|\bxi|^2 - d \ge 2\sqrt{d x} + 2x\bigr) \le e^{-x}, \quad x > 0.
    \end{equation}
    Let $\sigma \in \bigl(0, r / (8\sqrt{d})\bigr]$ and $x = r^2 / (16\sigma^2)$, yielding
    \begin{equation*}
        d + 2\sqrt{d x} + 2x \le \frac{r^2}{64\sigma^2} + \frac{r^2}{16\sigma^2} + \frac{r^2}{8\sigma^2} = \frac{13 r^2}{64 \sigma^2} \le \frac{r^2}{4 \sigma^2},
    \end{equation*}
    which further implies
    \begin{equation*}
        A_{\sigma}^c = \bigl\{|\bxi| > r / (2\sigma)\bigr\} \subseteq \bigl\{|\bxi|^2 \ge d + 2\sqrt{d x} + 2x\bigr\}.
    \end{equation*}
    Thus we have
    \begin{equation}\label{eq_PAc}
        \mathbb{P}(A_{\sigma}^c) \le \mathbb{P}\bigl(|\bxi|^2 - d \ge 2\sqrt{d x} + 2x\bigr) \le \exp\bigl(-r^2 / (16\sigma^2)\bigr) \leq \exp\bigl(-r^2 / (16\bar{\sigma}^2)\bigr)
    \end{equation}
    for every $\sigma \in (0, \bar{\sigma}]$, where the second inequality follows from \eqref{eq0_lm} and the definition of $x$.

    For later use, we introduce the following inequality:
    \begin{equation}\label{eq_mono}
        \sigma^{-s} \exp\bigl(-r^2 / (32\sigma^2)\bigr) \le \bar{\sigma}^{-s} \exp\bigl(-r^2 / (32\bar{\sigma}^2)\bigr), \quad \sigma \in (0, \bar{\sigma}], \quad s \in \{1, \ldots, 2p + n\},
    \end{equation}
    which follows from the fact that the function $\sigma \mapsto \sigma^{-s} \exp(-r^2 / (32\sigma^2))$ is increasing on $(0, \bar{\sigma}]$ for every $s \in \{1, \ldots, 2p + n\}$.

    \subsubsection*{(ii) Bound on the bias}

    To prove \eqref{eq_thm_mean}, we take expectations of both sides of \eqref{eq_app_dec} and insert the following identities
    \begin{align*}
        \E\bigl[\mathcal{E}^{\alpha}[u]\,\mathbf{1}_{A_{\sigma}}\bigr] & = \E\bigl[\mathcal{E}^{\alpha}[u]\bigr] - \E\bigl[\mathcal{E}^{\alpha}[u]\,\mathbf{1}_{A_{\sigma}^c}\bigr], \\
        \E\bigl[H_{\alpha}(\bxi)T_p(\bx, \bxi)\,\mathbf{1}_{A_{\sigma}}\bigr] & = \E\bigl[H_{\alpha}(\bxi)T_p(\bx, \bxi)\bigr] - \E\bigl[H_{\alpha}(\bxi)T_p(\bx, \bxi)\,\mathbf{1}_{A_{\sigma}^c}\bigr]
    \end{align*}
    to deduce
    \begin{align}
        &\E\bigl[\mathcal{E}^{\alpha}[u]\bigr] - \E\bigl[H_{\alpha}(\bxi)T_p(\bx, \bxi)\bigr] \notag\\
        = \;& \sigma^n\,\E\bigl[H_{\alpha}(\bxi)\mathcal{R}(\bx, \bxi, \sigma)\,\mathbf{1}_{A_{\sigma}}\bigr] + \E\bigl[\mathcal{E}^{\alpha}[u]\,\mathbf{1}_{A_{\sigma}^c}\bigr] - \E\bigl[H_{\alpha}(\bxi)T_p(\bx, \bxi)\,\mathbf{1}_{A_{\sigma}^c}\bigr] \label{eq_mean_dec}
    \end{align}
    for every $\sigma > 0$.
    Now we simplify \eqref{eq_mean_dec} to obtain \eqref{eq_thm_mean}.

    The second term $\E\bigl[H_{\alpha}(\bxi)T_p(\bx, \bxi)\bigr]$ in \eqref{eq_mean_dec} admits the following simplification:
    \begin{equation}\label{eq_mom_lead}
        \E\bigl[H_{\alpha}(\bxi) T_p(\bx, \bxi)\bigr] = \sum_{|\beta| = p} \frac{\partial^{\beta} u(\bx)}{\beta!}\,\E\bigl[H_{\alpha}(\bxi)\bxi^{\beta}\bigr] = \partial^{\alpha} u(\bx),
    \end{equation}
    where the first equality follows from \eqref{eq_taylor} and the second equality follows from \eqref{eq_ortho}.
    The first term on the right-hand side of \eqref{eq_mean_dec} can be bounded by
    \begin{equation}\label{eq_mean_rem}
        \begin{aligned}
            \bigl|\E\bigl[H_{\alpha}(\bxi)\mathcal{R}(\bx, \bxi, \sigma)\,\mathbf{1}_{A_{\sigma}}\bigr]\bigr| & \le \E\bigl[\bigl|H_{\alpha}(\bxi)\mathcal{R}(\bx, \bxi, \sigma)\bigr|\,\mathbf{1}_{A_{\sigma}}\bigr] \\
            & \le \E\bigl[\bigl|H_{\alpha}(\bxi)\mathcal{R}(\bx, \bxi, \sigma)\bigr|^q\,\mathbf{1}_{A_{\sigma}}\bigr]^{1 / q} \\
            & \le \Bigl[\bigl(\E\bigl[|H_{\alpha}(\bxi)|^{2q}\bigr]\bigr)^{1 / 2} \bigl(\E\bigl[|\mathcal{R}(\bx, \bxi, \sigma)|^{2q}\,\mathbf{1}_{A_{\sigma}}\bigr]\bigr)^{1 / 2}\Bigr]^{1 / q}
        \end{aligned}
    \end{equation}
    for every $q \ge 1$, where the second and the third inequalities follow from H\"older's inequality and the Cauchy--Schwarz inequality, respectively.

    To bound the second factor in the last line of \eqref{eq_mean_rem}, we note that for every $\bm{z} \in \overline{D_{r / 2}}$ and every nonzero $\bm{v} \in \mathbb{C}^d$, the function $s \mapsto u(\bm{z} + s\bm{v})$ is analytic on $\{s \in \mathbb{C} : |s| < r / (2|\bm{v}|)\}$ and bounded by $M$ in modulus, so the Cauchy estimate gives
    \begin{equation}\label{eq_cauchy}
        \bigl|(\bm{v} \cdot \nabla)^{n^{\prime}} u(\bm{z})\bigr| \le n^{\prime}!\, M (2 / r)^{n^{\prime}}\, |\bm{v}|^{n^{\prime}}, \quad \bm{z} \in \overline{D_{r / 2}}, \ \bm{v} \in \mathbb{C}^d.
    \end{equation}
    Inserting \eqref{eq_cauchy} into \eqref{eq_rem_int} and using $\int_0^1 (1 - \theta)^{n^{\prime} - 1} / (n^{\prime} - 1)!\, d\theta = 1 / n^{\prime}!$ gives
    \begin{equation}\label{eq_rem_pt}
        \bigl|\mathcal{R}(\bx, \bxi, \sigma)\bigr|\,\mathbf{1}_{A_{\sigma}} \le M (2 / r)^{n^{\prime}}\, |\bxi|^{n^{\prime}}, \quad \sigma > 0.
    \end{equation}
    Combining \eqref{eq_rem_pt} with the Cauchy--Schwarz inequality gives
    \begin{equation}\label{eq_rem_bd}
        \E\bigl[\bigl|H_{\alpha}(\bxi)\,\mathcal{R}(\bx, \bxi, \sigma)\bigr|^q\,\mathbf{1}_{A_{\sigma}}\bigr] \le M^q (2 / r)^{qn^{\prime}}\, \bigl(\E\bigl[|H_{\alpha}(\bxi)|^{2q}\bigr]\bigr)^{1 / 2}\, \E\bigl[|\bxi|^{2 qn^{\prime}}\bigr]^{1 / 2},
    \end{equation}
    for all $\sigma > 0$ and every $q \ge 1$.

    The second term on the right-hand side of \eqref{eq_mean_dec} can be bounded as follows:
    \begin{equation}\label{eq_bad_cs}
        \begin{aligned}
            \E\Bigl[\bigl| \mathcal{E}^{\alpha}[u] \bigr|^q \mathbf{1}_{A_{\sigma}^c} \Bigr] & \le (2M)^q \sigma^{-pq}\,\E\Bigl[\bigl| H_{\alpha}(\bxi) \bigr|^q\,\mathbf{1}_{A_{\sigma}^c} \Bigr] \\
            & \le (2M)^q \sigma^{-pq}\, \bigl(\E\bigl[|H_{\alpha}(\bxi)|^{2q}\bigr]\bigr)^{1 / 2} \, \mathbb{P}(A_{\sigma}^c)^{1 / 2} \\
            & \le \sigma^n (2M)^q\, \bigl(\E\bigl[|H_{\alpha}(\bxi)|^{2q}\bigr]\bigr)^{1 / 2}\, \bar{\sigma}^{-pq - n} \exp\bigl(-r^2 / (32\bar{\sigma}^2)\bigr)
        \end{aligned}
    \end{equation}
    for every $q \in \{1, 2\}$ and all $\sigma \in (0, \bar{\sigma})$, where the first inequality is obtained by inserting \eqref{eq_sing_stein} into \eqref{eq_rsse} and using the bound $|u| \le M$ assumed in Assumption~\ref{assu_main}; the second inequality follows from the Cauchy--Schwarz inequality; the last inequality follows from \eqref{eq_PAc} and \eqref{eq_mono} with $s = pq + n \le 2p + n$.
    Similarly, by the Cauchy--Schwarz inequality, \eqref{eq_PAc}, and \eqref{eq_mono} with $s = n$, the third term on the right-hand side of \eqref{eq_mean_dec} can be bounded as follows:
    \begin{equation}\label{eq_bdEpsbad}
        \begin{aligned}
            \E\bigl[\bigl|H_{\alpha}(\bxi)T_p(\bx, \bxi)\bigr|^q\,\mathbf{1}_{A_{\sigma}^c}\bigr] & \le \bigl(\E\bigl[|H_{\alpha}(\bxi)T_p(\bx, \bxi)|^{2q}\bigr]\bigr)^{1 / 2}\, \exp\bigl(-r^2 / (32\sigma^2)\bigr) \\
            & \le \sigma^n \bigl(\E\bigl[|H_{\alpha}(\bxi)T_p(\bx, \bxi)|^{2q}\bigr]\bigr)^{1 / 2}\, \bar{\sigma}^{-n} \exp\bigl(-r^2 / (32\bar{\sigma}^2)\bigr)
        \end{aligned}
    \end{equation}
    for every $q \in \{1, 2\}$ and all $\sigma \in (0, \bar{\sigma})$.

    Taking absolute values on both sides of \eqref{eq_mean_dec}, applying the triangle inequality, and inserting \eqref{eq_mom_lead} and the $q = 1$ cases of \eqref{eq_mean_rem}, \eqref{eq_rem_bd}, \eqref{eq_bad_cs}, and \eqref{eq_bdEpsbad} into the resulting inequality, we obtain \eqref{eq_thm_mean} with
    \begin{equation}\label{eq_cmean}
        \begin{aligned}
            C_{\mathrm{E}} := \;& M (2 / r)^{n^{\prime}}\, \bigl(\E\bigl[|H_{\alpha}(\bxi)|^2\bigr]\bigr)^{1 / 2}\, \bigl(\E\bigl[|\bxi|^{2n^{\prime}}\bigr]\bigr)^{1 / 2} \\
            & + 2M\, \bar{\sigma}^{-p - n}\, \bigl(\E\bigl[|H_{\alpha}(\bxi)|^2\bigr]\bigr)^{1 / 2}\, \exp\bigl(-r^2 / (32\bar{\sigma}^2)\bigr) \\
            & + \bar{\sigma}^{-n}\, \bigl(\E\bigl[\bigl|H_{\alpha}(\bxi)T_p(\bx, \bxi)\bigr|^2\bigr]\bigr)^{1 / 2}\, \exp\bigl(-r^2 / (32\bar{\sigma}^2)\bigr) \\
            \le \;& \frac{3 M (\alpha!)^{1 / 2}}{p!}\, \bigl(4 (d + n^{\prime})\bigr)^{2 n^{\prime}}\, r^{-n^{\prime}}.
        \end{aligned}
    \end{equation}

    \subsubsection*{(iii) Bound on the variance}

    To prove \eqref{eq_thm_var}, we use the identity $\Var[Z] = \E[|Z|^2] - |\E[Z]|^2$, valid for any complex-valued random variable $Z$, to obtain
    \begin{equation}\label{eq_varvar}
        \begin{aligned}
            &\Var\bigl[\mathcal{E}^{\alpha}[u](\bx; \bxi)\bigr] - \Var\bigl[H_{\alpha}(\bxi) \sum_{|\beta| = p} \frac{\bxi^{\beta}}{\beta!}\,\partial^{\beta} u(\bx) \bigr] \\
            = \;& \E\bigl[|\mathcal{E}^{\alpha}[u]|^2\bigr] - \E\bigl[|H_{\alpha}(\bxi)T_p(\bx, \bxi)|^2\bigr] - \Bigl| \E\bigl[\mathcal{E}^{\alpha}[u](\bx; \bxi) \bigr] \Bigr|^2 + \Bigl| \E\bigl[H_{\alpha}(\bxi)T_p(\bx, \bxi)\bigr] \Bigr|^2,
        \end{aligned}
    \end{equation}
    where $T_p(\bx, \bxi)$ is defined by \eqref{eq_taylor}.

    To bound the first two terms on the right-hand side, we take $\E[| \cdot |^2 \mathbf{1}_{A_{\sigma}}]$ on both sides of \eqref{eq_app_dec} to obtain
    \begin{equation*}
        \begin{aligned}
            \E\bigl[|\mathcal{E}^{\alpha}[u]|^2\,\mathbf{1}_{A_{\sigma}}\bigr] = \;& \E\bigl[|H_{\alpha}(\bxi)T_p(\bx, \bxi)|^2\,\mathbf{1}_{A_{\sigma}}\bigr] + \sigma^{2n}\,\E\bigl[\bigl|H_{\alpha}(\bxi)\mathcal{R}(\bx, \bxi, \sigma)\bigr|^2\,\mathbf{1}_{A_{\sigma}}\bigr]\\
            & + 2\sigma^n\,\mathrm{Re}\Big\{\E\bigl[H_{\alpha}(\bxi)T_p(\bx, \bxi)\,\overline{H_{\alpha}(\bxi)\mathcal{R}(\bx, \bxi, \sigma)}\,\mathbf{1}_{A_{\sigma}}\bigr] \Big\}.
        \end{aligned}
    \end{equation*}
    Splitting the expectations over $A_{\sigma}$ and $A_{\sigma}^c$, the above identity is equivalent to
    \begin{align}
        &\E\bigl[|\mathcal{E}^{\alpha}[u]|^2\bigr] - \E\bigl[|H_{\alpha}(\bxi)T_p(\bx, \bxi)|^2\bigr] \notag\\
        = \;& 2\sigma^n\,\mathrm{Re}\Big\{\E\bigl[H_{\alpha}(\bxi)T_p(\bx, \bxi)\,\overline{H_{\alpha}(\bxi)\mathcal{R}(\bx, \bxi, \sigma)}\,\mathbf{1}_{A_{\sigma}}\bigr] \Big\} + \sigma^{2n}\,\E\bigl[\bigl|H_{\alpha}(\bxi)\mathcal{R}(\bx, \bxi, \sigma)\bigr|^2\,\mathbf{1}_{A_{\sigma}}\bigr] \notag\\
        & + \E\bigl[|\mathcal{E}^{\alpha}[u]|^2\,\mathbf{1}_{A_{\sigma}^c}\bigr] - \E\bigl[|H_{\alpha}(\bxi)T_p(\bx, \bxi)|^2\,\mathbf{1}_{A_{\sigma}^c}\bigr] \label{eq_mom2_dec}
    \end{align}
    for all $\sigma \in (0, \bar{\sigma})$.
    For rigor, we remark that each expectation in \eqref{eq_mom2_dec} is finite by the $q = 2$ cases of \eqref{eq_rem_bd}, \eqref{eq_bad_cs}, and \eqref{eq_bdEpsbad} and the Cauchy--Schwarz inequality.
    Now we bound the terms on the right-hand side of \eqref{eq_mom2_dec}.

    The first term on the right-hand side of \eqref{eq_mom2_dec} can be bounded by
    \begin{equation}\label{eq0_cross_bd}
        \begin{aligned}
            &2\sigma^n\,\Bigl|\E\bigl[H_{\alpha}(\bxi)T_p(\bx, \bxi)\,\overline{H_{\alpha}(\bxi)\mathcal{R}(\bx, \bxi, \sigma)}\,\mathbf{1}_{A_{\sigma}}\bigr]\Bigr| \\
            \le \;& 2\sigma^n\,\bigl(\E\bigl[\bigl|H_{\alpha}(\bxi)T_p(\bx, \bxi)\bigr|^2\bigr]\bigr)^{1 / 2}\,\bigl(\E\bigl[\bigl|H_{\alpha}(\bxi)\mathcal{R}(\bx, \bxi, \sigma)\bigr|^2\,\mathbf{1}_{A_{\sigma}}\bigr]\bigr)^{1 / 2} \\
            \le \;& 2M (2 / r)^{n^{\prime}}\,\bigl(\E\bigl[\bigl|H_{\alpha}(\bxi)T_p(\bx, \bxi)\bigr|^2\bigr]\bigr)^{1 / 2}\,\bigl(\E\bigl[|H_{\alpha}(\bxi)|^4\bigr]\bigr)^{1 / 4}\,\bigl(\E\bigl[|\bxi|^{4n^{\prime}}\bigr]\bigr)^{1 / 4}\, \sigma^n
        \end{aligned}
    \end{equation}
    for all $\sigma \in (0, \bar{\sigma})$, where the first inequality follows from the Cauchy--Schwarz inequality and the second inequality follows from \eqref{eq_rem_bd} with $q = 2$.

    The second term on the right-hand side of \eqref{eq_mom2_dec} can be bounded by
    \begin{equation}\label{eq0_rem_sq}
        \sigma^{2n}\,\E\bigl[\bigl|H_{\alpha}(\bxi)\mathcal{R}(\bx, \bxi, \sigma)\bigr|^2\,\mathbf{1}_{A_{\sigma}}\bigr] \le M^2 (2 / r)^{2n^{\prime}}\,\bigl(\E\bigl[|H_{\alpha}(\bxi)|^4\bigr]\bigr)^{1 / 2}\,\bigl(\E\bigl[|\bxi|^{4n^{\prime}}\bigr]\bigr)^{1 / 2}\, \bar{\sigma}^n\, \sigma^n
    \end{equation}
    for all $\sigma \in (0, \bar{\sigma})$, which follows from \eqref{eq_rem_bd} with $q = 2$ and the inequality $\sigma^{2n} \le \bar{\sigma}^n\, \sigma^n$.
    The third and fourth terms on the right-hand side of \eqref{eq_mom2_dec} are bounded by the $q = 2$ cases of \eqref{eq_bad_cs} and \eqref{eq_bdEpsbad}, respectively, which gives
    \begin{align}
        \E\bigl[\bigl|\mathcal{E}^{\alpha}[u]\bigr|^2\,\mathbf{1}_{A_{\sigma}^c}\bigr] & \le 4M^2\,\bigl(\E\bigl[|H_{\alpha}(\bxi)|^4\bigr]\bigr)^{1 / 2}\, \bar{\sigma}^{-2p - n}\, \exp\bigl(-r^2 / (32\bar{\sigma}^2)\bigr)\, \sigma^n, \label{eq0_tail_sq}\\
        \E\bigl[\bigl|H_{\alpha}(\bxi)T_p(\bx, \bxi)\bigr|^2\,\mathbf{1}_{A_{\sigma}^c}\bigr] & \le \bigl(\E\bigl[\bigl|H_{\alpha}(\bxi)T_p(\bx, \bxi)\bigr|^4\bigr]\bigr)^{1 / 2}\, \bar{\sigma}^{-n}\, \exp\bigl(-r^2 / (32\bar{\sigma}^2)\bigr)\, \sigma^n \label{eq1_tail_sq}
    \end{align}
    for all $\sigma \in (0, \bar{\sigma})$.
    Taking $| \cdot |$ on both sides of \eqref{eq_mom2_dec}, applying the triangle inequality, and inserting \eqref{eq0_cross_bd}, \eqref{eq0_rem_sq}, and \eqref{eq0_tail_sq}--\eqref{eq1_tail_sq} into the resulting inequality, we obtain
    \begin{equation}\label{eq_cmom}
        \Bigl|\E\bigl[|\mathcal{E}^{\alpha}[u]|^2\bigr] - \E\bigl[|H_{\alpha}(\bxi)T_p(\bx, \bxi)|^2\bigr]\Bigr| \le C_{\mathrm{M}}\, \sigma^n
    \end{equation}
    for all $\sigma \in (0, \bar{\sigma})$, with
    \begin{equation}\label{eq_cmom_const}
        \begin{aligned}
            C_{\mathrm{M}} := \;& 2M (2 / r)^{n^{\prime}}\, \bigl(\E\bigl[\bigl|H_{\alpha}(\bxi)T_p(\bx, \bxi)\bigr|^2\bigr]\bigr)^{1 / 2}\, \bigl(\E\bigl[|H_{\alpha}(\bxi)|^4\bigr]\bigr)^{1 / 4}\, \bigl(\E\bigl[|\bxi|^{4n^{\prime}}\bigr]\bigr)^{1 / 4} \\
            & + M^2 (2 / r)^{2n^{\prime}}\, \bar{\sigma}^n\, \bigl(\E\bigl[|H_{\alpha}(\bxi)|^4\bigr]\bigr)^{1 / 2}\, \bigl(\E\bigl[|\bxi|^{4n^{\prime}}\bigr]\bigr)^{1 / 2} \\
            & + 4M^2\, \bar{\sigma}^{-2p - n}\, \bigl(\E\bigl[|H_{\alpha}(\bxi)|^4\bigr]\bigr)^{1 / 2}\, \exp\bigl(-r^2 / (32\bar{\sigma}^2)\bigr) \\
            & + \bar{\sigma}^{-n}\, \bigl(\E\bigl[\bigl|H_{\alpha}(\bxi)T_p(\bx, \bxi)\bigr|^4\bigr]\bigr)^{1 / 2}\, \exp\bigl(-r^2 / (32\bar{\sigma}^2)\bigr) \\
            \le \;& \frac{8 M^2 \alpha!}{(p!)^2}\, \bigl(3456 (d + n^{\prime})^4\bigr)^{n^{\prime}} r^{-(2p + n)}.
        \end{aligned}
    \end{equation}

    Using \eqref{eq_thm_mean} and \eqref{eq_mom_lead}, the last two terms on the right-hand side of \eqref{eq_varvar} can be bounded by
    \begin{equation*}
        \Bigl| \Bigl| \E\bigl[\mathcal{E}^{\alpha}[u](\bx; \bxi)\bigr] \Bigr|^2 - \Bigl| \E\bigl[H_{\alpha}(\bxi)T_p(\bx, \bxi)\bigr] \Bigr|^2 \Bigr| \leq 2\,|\partial^{\alpha} u(\bx)|\,C_{\mathrm{E}}\, \sigma^n + C_{\mathrm{E}}^2\, \sigma^{2n}.
    \end{equation*}
    Inserting this inequality and \eqref{eq_cmom} into \eqref{eq_varvar} gives \eqref{eq_thm_var} with
    \begin{equation}\label{eq_cvar}
        C_{\mathrm{V}} := C_{\mathrm{M}} + 2\,|\partial^{\alpha} u(\bx)|\,C_{\mathrm{E}} + \bar{\sigma}^n C_{\mathrm{E}}^2 \le \; \frac{23 M^2 \alpha!}{(p!)^2}\, \bigl(3456 (d + n^{\prime})^4\bigr)^{n^{\prime}} r^{-(2p + n)}.
    \end{equation}
\end{proof}


\section{Numerical experiments}\label{sec_experiments}

We first validate the main results in Theorem~\ref{thm_main} by examining the mean and variance of the RSSE derivative approximations, and then compare RSSE with antithetic sampling and automatic differentiation in PINNs.

\subsection{Derivative approximation}\label{sec_mom_val}


We consider the derivative approximation of the following test function:
\begin{equation}\label{eq_test_fn}
    u(\bx) = \frac{1}{\sqrt{3}}\sum_{i = 1}^{3} \sin\!\left(x_i + \cos x_{i + 1}\right), \quad \bx = (x_1, x_2, x_3), \quad x_4 := x_1.
\end{equation}
We evaluate its derivatives at $\bx_0 = (1 / 4, 1 / 4, 1 / 4)$ under the parameter settings $\alpha \in \{(1, 1, 0), (3, 0, 0), (2, 2, 0)\}$, $n \in \{4, 6, 8\}$, and $\sigma = 0.05 \times 2^{j / 5}$, $j = 0, \ldots, 5$.
For each $(\alpha, n, \sigma)$, we use $R = 25600$ independent seeds, with $K = 2^{15}$ samples per seed in \eqref{eq_ksamp_rsse}.
The estimate from seed $s$ is
\begin{equation}\label{eq_grp_est}
    Y_s^{\alpha, n}(\sigma) := \frac{1}{K}\sum_{\ell = 1}^{K} \mathcal{E}^{\alpha}[u](\bx_0; \bxi_{s, \ell}), \quad \bxi_{s, \ell}\stackrel{\mathrm{i.i.d.}}{\sim}\N(0, I_3).
\end{equation}
The resulting sample mean and the estimated variance of one RSSE evaluation are
\begin{equation}\label{eq_mom_est}
    \widehat{\mu}_R^{\alpha, n}(\sigma) := \frac{1}{R}\sum_{s = 1}^{R}Y_s^{\alpha, n}(\sigma), \quad \widehat{V}_R^{\alpha, n}(\sigma) := \frac{K}{R - 1}\sum_{s = 1}^{R} \left|Y_s^{\alpha, n}(\sigma) - \widehat{\mu}_R^{\alpha, n}(\sigma)\right|^2.
\end{equation}
Here, the factor $K$ rescales the variance of each sample average to the variance of one RSSE evaluation.

For each $(\alpha, n)$, we fit
\begin{equation}\label{eq_asy_fit}
    \widehat{\mu}_R^{\alpha, n}(\sigma) = a_{\mu} + b_{\mu} \sigma^{q_{\mu}}, \quad \widehat{V}_R^{\alpha, n}(\sigma) = a_V + b_V \sigma^{q_V}.
\end{equation}
We estimate these parameters by nonlinear least squares.
The notation $V_{0, \alpha}$ denotes the leading variance term of the RSSE, i.e., the second variance term inside $\lvert \cdot \rvert$ on the left-hand side of \eqref{eq_thm_var}, evaluated at $\bx = \bx_0$.

\begin{figure}[ht]
    \centering
    \newcommand{\rssemomentplot}[2]{%
        \begin{minipage}[t]{0.325\linewidth}
            \centering
            {\fontsize{5}{6}\selectfont\strut #2\par}
            \vspace{0.6pt}%
            \includegraphics[width=\linewidth,trim=0 0 0 17.28bp,clip]{#1}%
        \end{minipage}%
    }
    \rssemomentplot{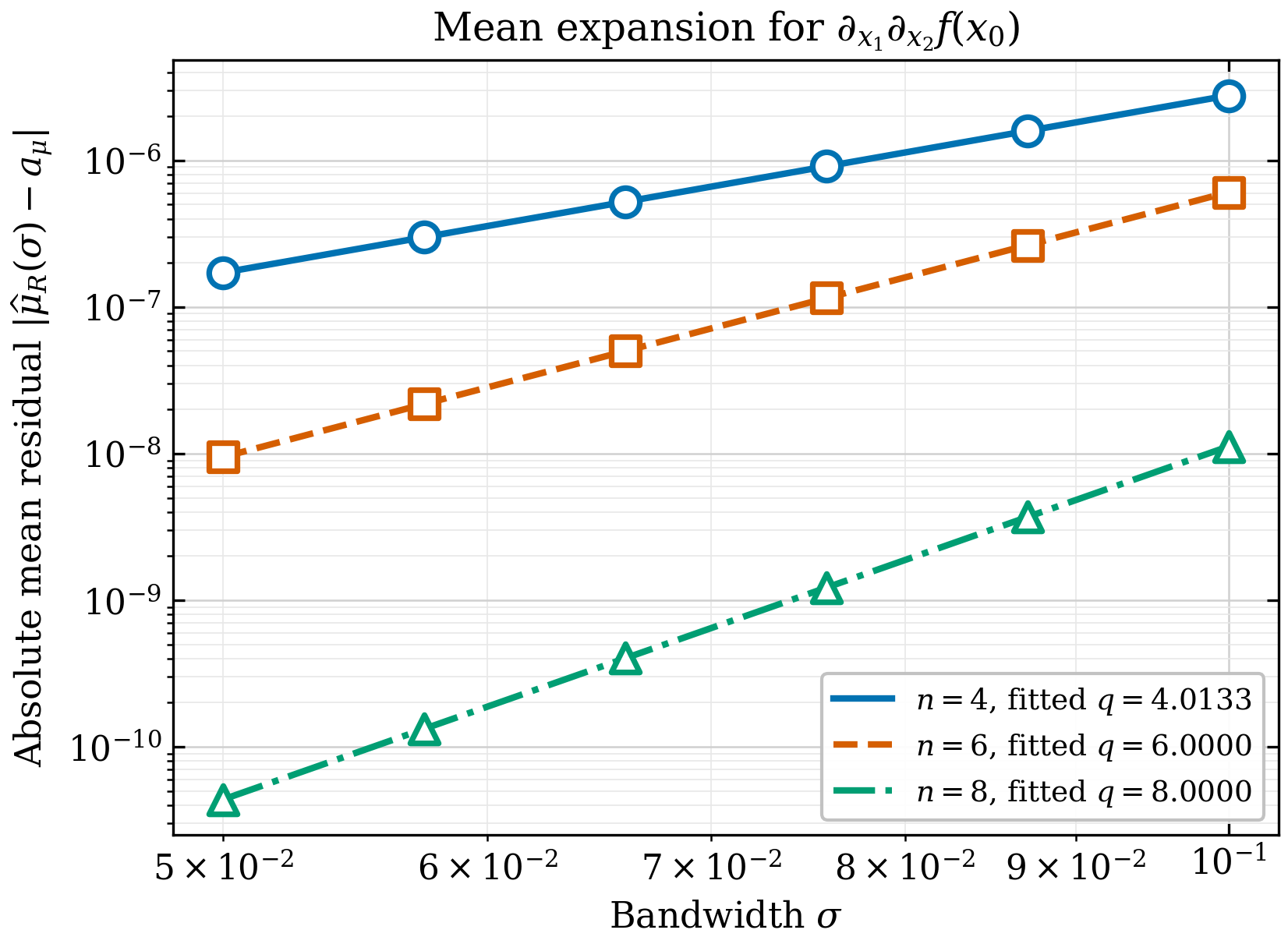}{Mean expansion for $\partial_{x_1}\partial_{x_2}u(\bx_0)$}\hfill
    \rssemomentplot{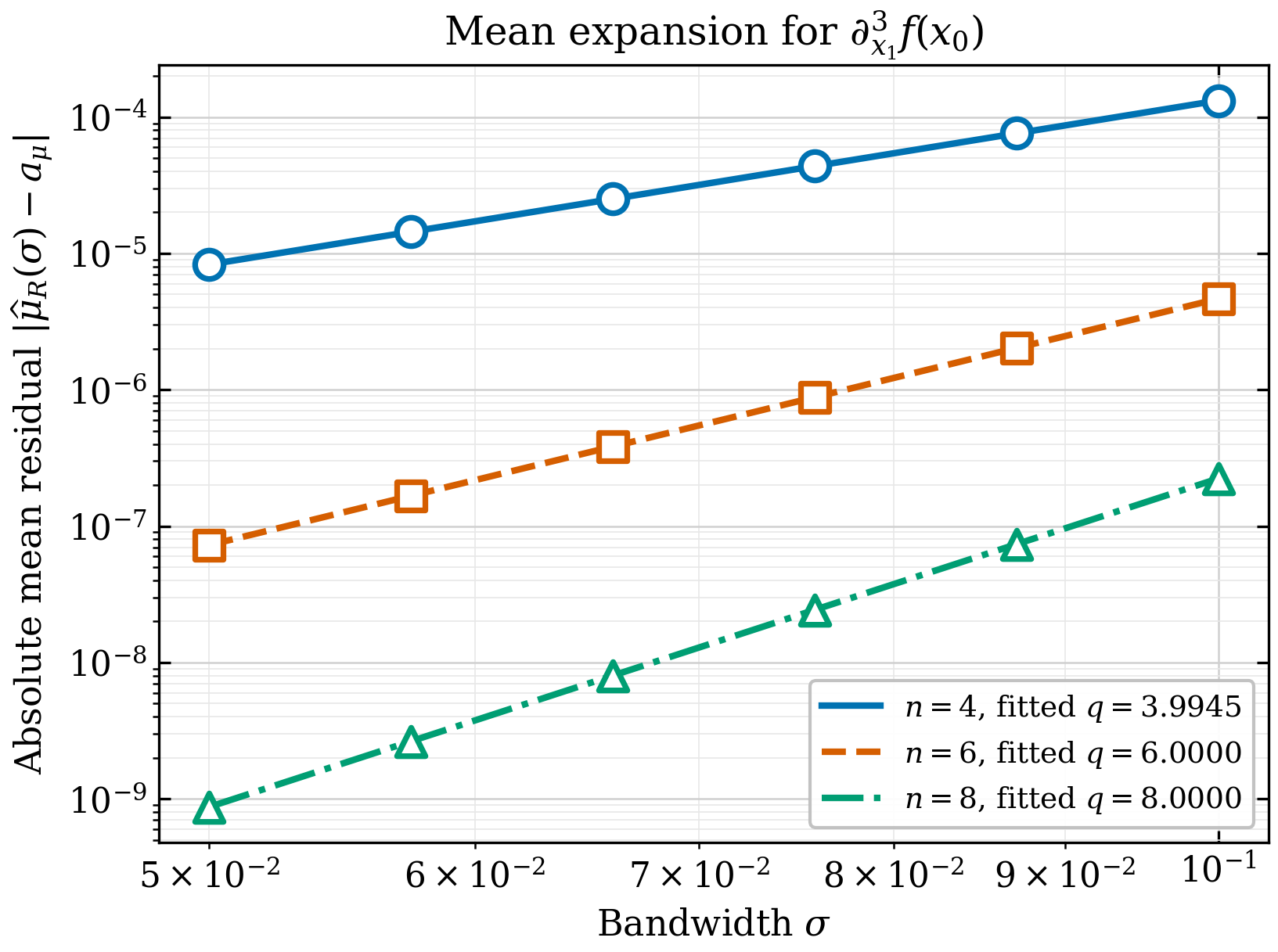}{Mean expansion for $\partial_{x_1}^3u(\bx_0)$}\hfill
    \rssemomentplot{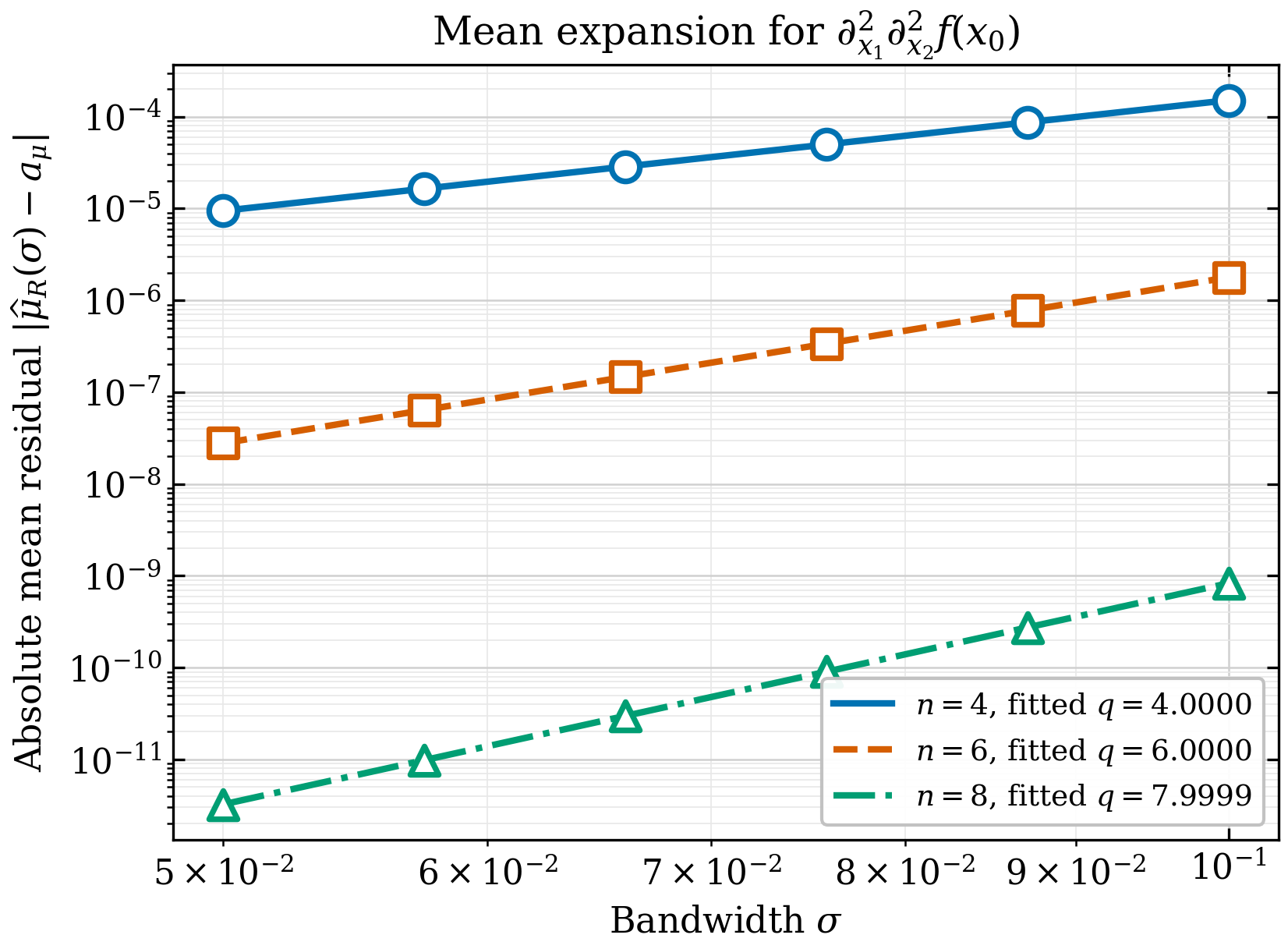}{Mean expansion for $\partial_{x_1}^2\partial_{x_2}^2u(\bx_0)$}\par\smallskip
    \rssemomentplot{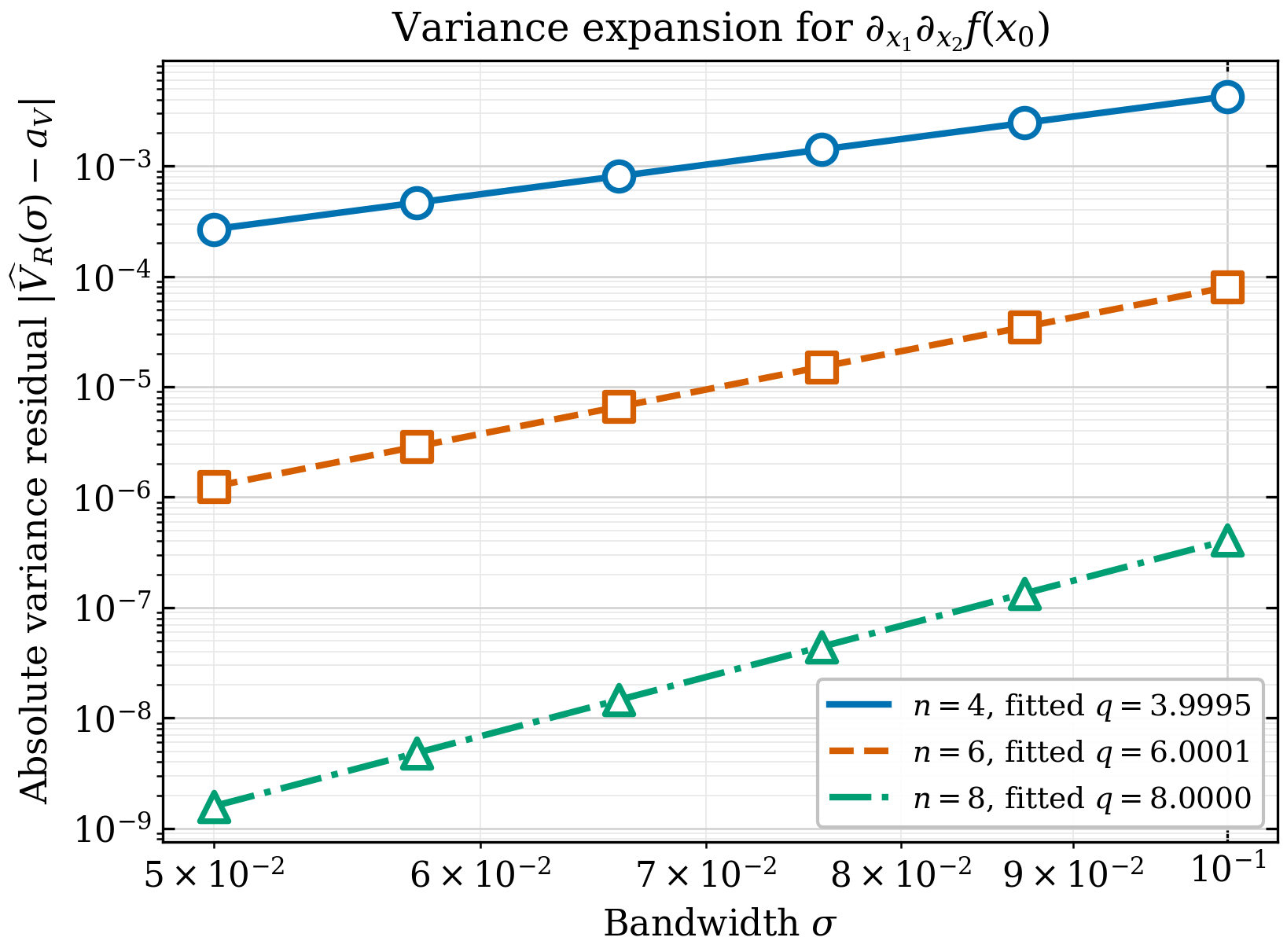}{Variance expansion for $\partial_{x_1}\partial_{x_2}u(\bx_0)$}\hfill
    \rssemomentplot{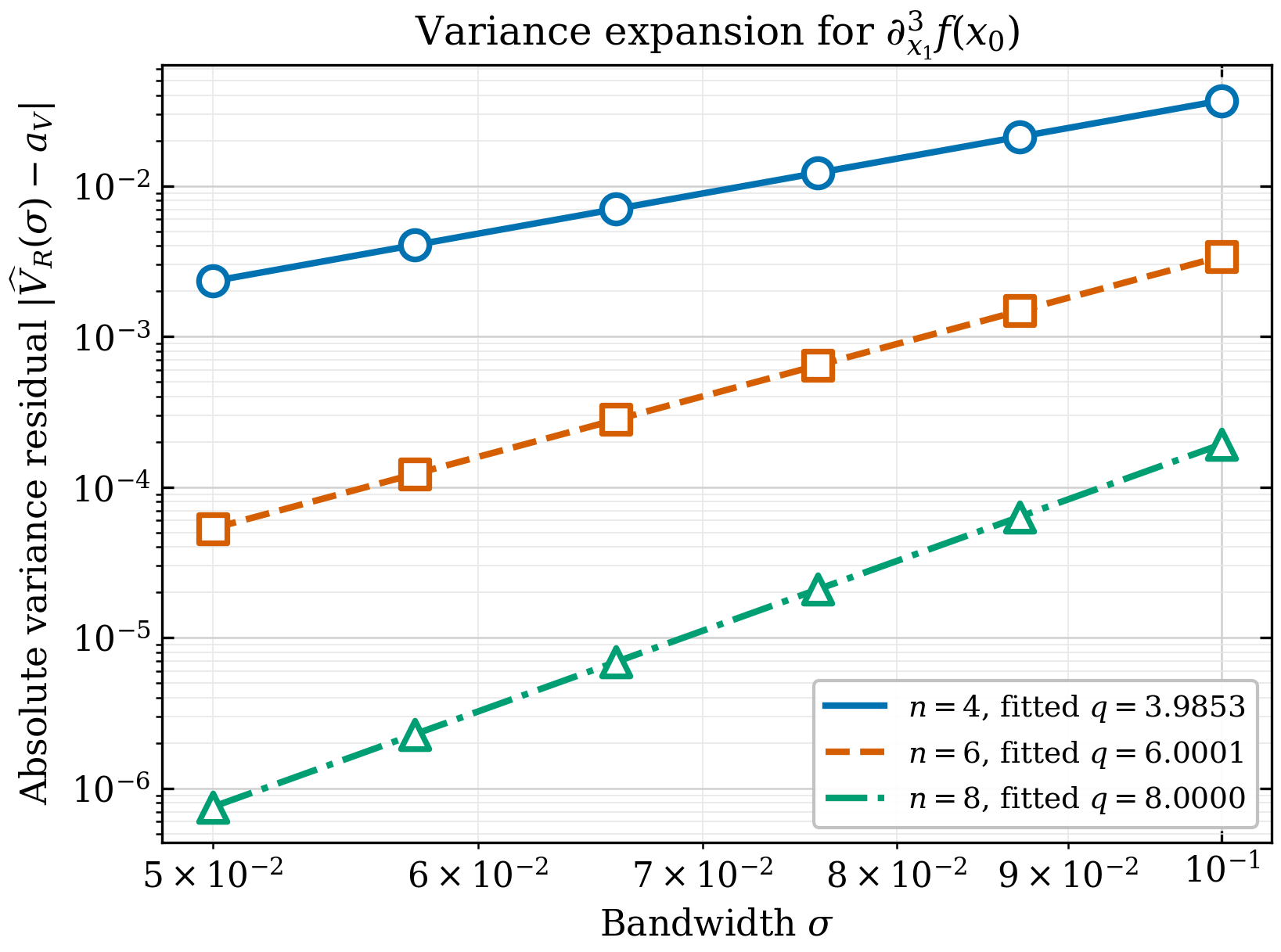}{Variance expansion for $\partial_{x_1}^3u(\bx_0)$}\hfill
    \rssemomentplot{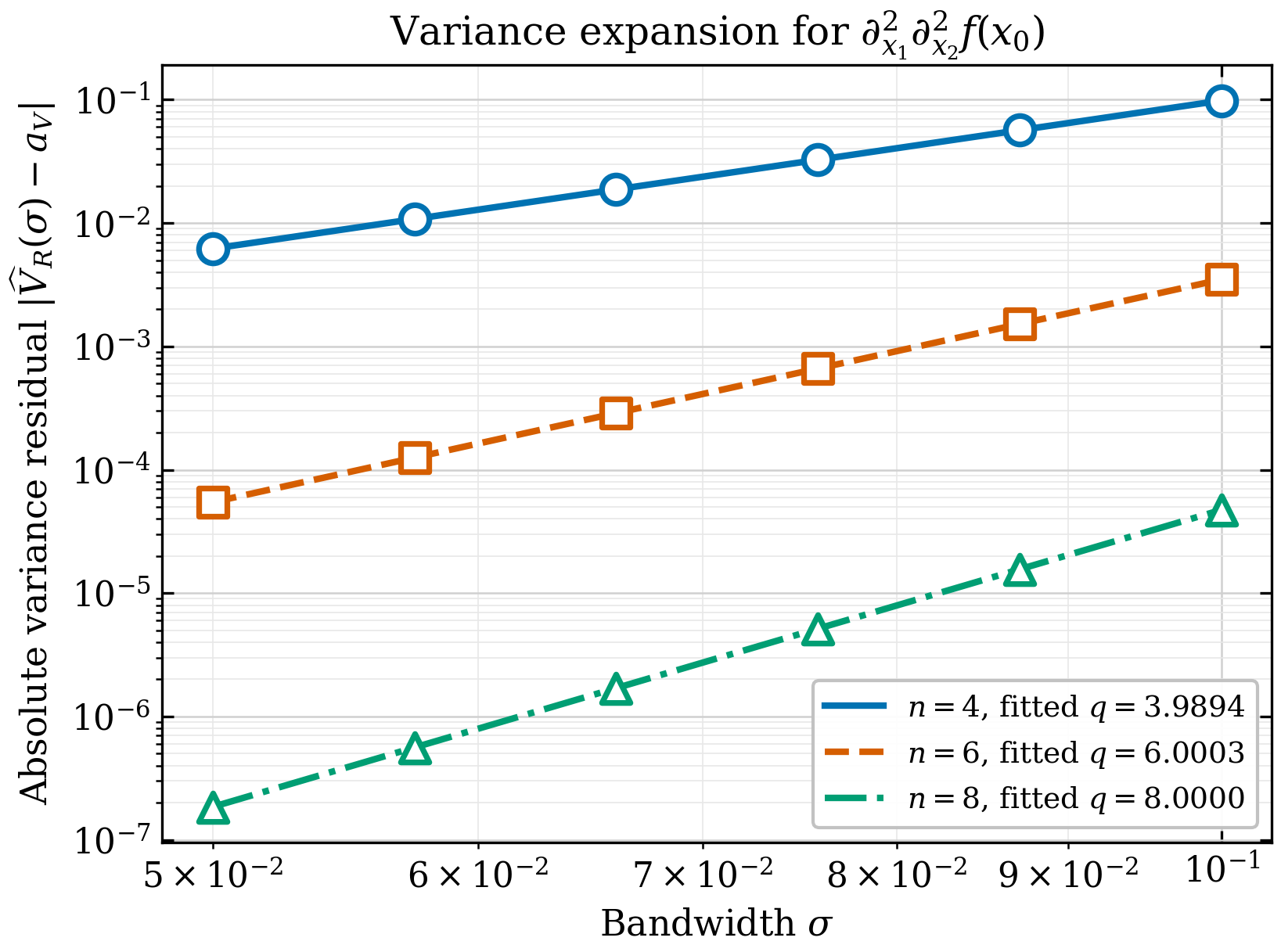}{Variance expansion for $\partial_{x_1}^2\partial_{x_2}^2u(\bx_0)$}
    \caption{Validation of the mean and variance expansions~\eqref{eq_thm_mean} and~\eqref{eq_thm_var} for the test function~\eqref{eq_test_fn}.
        The absolute mean residual is shown in the top row and the absolute variance residual in the bottom row for $\alpha = (1, 1, 0)$, $(3, 0, 0)$, and $(2, 2, 0)$ (left to right).
        Markers show the six step sizes for $n = 4, 6, 8$; lines show $\lvert b_{\mu}\rvert \sigma^{q_{\mu}}$ in the top row and $\lvert b_V\rvert \sigma^{q_V}$ in the bottom row.}
    \label{fig_mom_val}
\end{figure}

Across all tested cases, the fitted mean and variance intercepts agree with the target derivatives and leading variance values to within $0.25\%$ and $0.81\%$, respectively; see Table~\ref{tab_mom_val}.
The fitted exponents for the mean and variance corrections are close to $n$, consistent with the predicted $\sigma^n$ scaling over the sampled step-size range.
Figure~\ref{fig_mom_val} displays the residuals relative to the fitted intercepts together with the corresponding fitted correction terms.

\begin{table}[ht]
    \centering
    \scriptsize
    \setlength{\tabcolsep}{2.5pt}
    \caption{Mean- and variance-expansion fits for the test function~\eqref{eq_test_fn}.
        Relative-error columns compare $a_{\mu}$ with $\partial^{\alpha} u(\bx_0)$ and $a_V$ with $V_{0, \alpha}$, respectively.}
    \label{tab_mom_val}
    \begin{tabular}{ccrrrrrrrr}
        \toprule
        $\alpha$ & $n$
        & \multicolumn{4}{c}{mean}
        & \multicolumn{4}{c}{variance} \\
        \cmidrule(lr){3-6}\cmidrule(lr){7-10}
        & & \multicolumn{1}{c}{$\partial^{\alpha} u(\bx_0)$}
        & \multicolumn{1}{c}{$a_{\mu}$}
        & \multicolumn{1}{c}{error (\%)}
        & \multicolumn{1}{c}{$q_{\mu}$}
        & \multicolumn{1}{c}{$V_{0, \alpha}$}
        & \multicolumn{1}{c}{$a_V$}
        & \multicolumn{1}{c}{error (\%)}
        & \multicolumn{1}{c}{$q_V$} \\
        \midrule
        $(1, 1, 0)$ & 4 &  0.1341 &  0.1340 & 0.05164 & 4.013 &  9.540 &  9.565 & 0.2589 & 4.000 \\
        $(1, 1, 0)$ & 6 &         &  0.1342 & 0.08839 & 6.000 &        &  9.578 & 0.3906 & 6.000 \\
        $(1, 1, 0)$ & 8 &         &  0.1342 & 0.06629 & 8.000 &        &  9.517 & 0.2416 & 8.000 \\
        \addlinespace[2pt]
        $(3, 0, 0)$ & 4 & -0.5365 & -0.5365 & 0.0008120 & 3.994 & 57.65 & 57.62 & 0.04687 & 3.985 \\
        $(3, 0, 0)$ & 6 &         & -0.5364 & 0.02758 & 6.000 &       & 57.18 & 0.8060 & 6.000 \\
        $(3, 0, 0)$ & 8 &         & -0.5368 & 0.05070 & 8.000 &       & 57.37 & 0.4830 & 8.000 \\
        \addlinespace[2pt]
        $(2, 2, 0)$ & 4 &  0.2260 &  0.2265 & 0.2455 & 4.000 & 47.55 & 47.23 & 0.6662 & 3.989 \\
        $(2, 2, 0)$ & 6 &         &  0.2262 & 0.07672 & 6.000 &       & 47.67 & 0.2538 & 6.000 \\
        $(2, 2, 0)$ & 8 &         &  0.2260 & 0.01167 & 8.000 &       & 47.48 & 0.1362 & 8.000 \\
        \bottomrule
    \end{tabular}
\end{table}

\FloatBarrier


\subsection{PDE solution}\label{sec_pinn_exp}

We compare the RSSE-PINN method in Section~\ref{sec_monte_carlo} with PINNs using antithetic sampling or automatic differentiation (AutoDiff) for polyharmonic and Cahn--Hilliard equations.
Both PDE studies use a fully connected network with four hidden layers, width 128, and sine activation, together with a batch size of 128.
Training uses AdamW~\cite{KingmaBa2015Adam,Loshchilov2019AdamW} for 1200 seconds with learning rate $10^{-4}$, weight decay $10^{-4}$, and gradient clipping at 1.

For RSSE and the antithetic estimator, we set $K = 256$ and test $\sigma \in \{0.05, 0.10, \ldots, 0.60\}$.
AutoDiff is independent of $\sigma$.
Each stochastic configuration is repeated five times with independent random seeds.
For the polyharmonic experiments in Section~\ref{sec_poly}, RSSE uses $n = 2k$, matching the PDE order; for the Cahn--Hilliard experiments in Section~\ref{sec_ch}, it uses $n = 2m$, matching the highest spatial derivative order.

For a fixed validation set $D_{\mathrm{val}}$, with $\mathbf{z}$ denoting a spatial or space--time point as appropriate, we report
\begin{equation}\label{eq_rel_err}
    \operatorname{RE} := \left( \sum_{\mathbf{z} \in D_{\mathrm{val}}}|u_{\theta}(\mathbf{z}) - u_*(\mathbf{z})|^2 \right)^{1 / 2} / \left( \sum_{\mathbf{z} \in D_{\mathrm{val}}}|u_*(\mathbf{z})|^2 \right)^{1 / 2}.
\end{equation}
For the polyharmonic experiments in Section~\ref{sec_poly}, $D_{\mathrm{val}}$ consists of independent samples $\bz = r\bm{v} / |\bm{v}|$, where $r \sim \operatorname{Unif}[0, 1]$ and $\bm{v} \sim \mathcal{N}(\mathbf{0}, I_d)$ are independent.
For the Cahn--Hilliard experiments in Section~\ref{sec_ch}, $D_{\mathrm{val}}$ consists of independent space--time samples $\bz = (\bx, t) \sim \operatorname{Unif}([0, 1]^d \times [0, 1])$.
In each case, the interior collocation samples in $D_{\mathrm{train}}$ and the validation samples in $D_{\mathrm{val}}$ are drawn independently and identically from the corresponding distribution.

For each method, problem setting, and seed, we take $\operatorname{RE}$ at the final recorded training time.
We report its mean and sample standard deviation across seeds.
For each stochastic method, the selected $\sigma$ minimizes the five-seed mean final $\operatorname{RE}$ over the tested candidates.

\subsubsection{Polyharmonic equations}\label{sec_poly}

For $k \in \{2, 3, 4\}$ and $d \in \{5, 10, 20\}$, let $\Omega := \{\bx \in \R^d: |\bx| < 1\}$ and $\Gamma := \partial\Omega = \{\bx \in \R^d: |\bx| = 1\}$.
We solve
\begin{equation}\label{eq_poly}
    \begin{cases}
        \Delta^k u(\bx) = f(\bx), & \bx \in \Omega, \\
        \partial_{\nu}^j u(\bx) = 0, & \bx \in \Gamma, \quad j = 0, \ldots, k - 1,
    \end{cases}
\end{equation}
where $\nu$ denotes the outward unit normal.
The notation $\partial_{\nu}^j u$ denotes the normal derivative of $u$ of order $j$.
We adopt the following hard-constraint structure to satisfy the zero boundary conditions automatically~\cite{Lu2021HardConstraints}:
\begin{equation}\label{eq_hardbc}
    u_{\theta}(\bx) := (1 - |\bx|^2)^k \psi_{\theta}(\bx),
\end{equation}
where $\psi_{\theta}(\bx)$ is the fully connected neural network and $u_{\theta}(\bx)$ is the resulting boundary-constrained approximation to $u_*(\bx)$.
The exact solution is
\begin{equation}\label{eq_exact}
    u_*(\bx) := (1 - |\bx|^2)^k \sum_{i = 1}^d\left(x_i + x_i^3\right)
\end{equation}
and $f := \Delta^k u_*$.

Table~\ref{tab_summary} reports the selected errors and step sizes, and Figures~\ref{fig_poly_k2}--\ref{fig_poly_k4} show the complete $\sigma$ sweeps and the corresponding training trajectories.
RSSE attains the smallest mean final $\operatorname{RE}$ for all nine $(k, d)$ pairs.
As $d$ increases, the antithetic estimator learns more slowly and sometimes deteriorates, whereas RSSE continues to reduce the error.
AutoDiff is less accurate in every completed setting and runs out of memory for $(k, d) = (4, 20)$.

\begingroup
    \addtolength{\intextsep}{-4pt}
    \addtolength{\floatsep}{-4pt}

    \begin{table}[htbp]
        \centering
        \scriptsize
        \setlength{\tabcolsep}{2.5pt}
        \caption{Numerical results for the polyharmonic PDE defined by \eqref{eq_poly}.
            For RSSE and antithetic sampling, we report the five-seed mean $ \pm $ sample standard deviation of the final $\operatorname{RE}$ defined in~\eqref{eq_rel_err}; $\sigma$ is selected to minimize the mean final $\operatorname{RE}$.
            The AutoDiff column reports the mean $ \pm $ sample standard deviation of the final $\operatorname{RE}$; \textsc{OOM} denotes an AutoDiff out-of-memory failure.}
        \label{tab_summary}
        \begin{tabular}{ccccccc}
            \toprule
            $k$ & $d$
            & \multicolumn{2}{c}{RSSE}
            & \multicolumn{2}{c}{antithetic}
            & \multicolumn{1}{c}{AutoDiff} \\
            \cmidrule(lr){3-4}\cmidrule(lr){5-6}\cmidrule(lr){7-7}
            & & \multicolumn{1}{c}{final $\operatorname{RE}$}
            & \multicolumn{1}{c}{$\sigma$}
            & \multicolumn{1}{c}{final $\operatorname{RE}$}
            & \multicolumn{1}{c}{$\sigma$}
            & \multicolumn{1}{c}{final $\operatorname{RE}$} \\
            \midrule
            2 & 5  & $5.765\mathrm{E}{-3} \pm 1.935\mathrm{E}{-3}$ & 0.15 & $2.553\mathrm{E}{-2} \pm 6.802\mathrm{E}{-3}$ & 0.15 & $9.579\mathrm{E}{-2} \pm 3.477\mathrm{E}{-3}$ \\
            2 & 10 & $4.285\mathrm{E}{-3} \pm 1.636\mathrm{E}{-3}$ & 0.15 & $5.403\mathrm{E}{-2} \pm 6.466\mathrm{E}{-3}$ & 0.20 & $4.391\mathrm{E}{-2} \pm 1.695\mathrm{E}{-3}$ \\
            2 & 20 & $1.178\mathrm{E}{-2} \pm 1.474\mathrm{E}{-3}$ & 0.20 & $3.715\mathrm{E}{-2} \pm 1.609\mathrm{E}{-3}$ & 0.15 & $4.876\mathrm{E}{-2} \pm 9.622\mathrm{E}{-4}$ \\
            \addlinespace[2pt]
            3 & 5  & $1.225\mathrm{E}{-2} \pm 2.340\mathrm{E}{-3}$ & 0.20 & $2.765\mathrm{E}{-2} \pm 6.025\mathrm{E}{-3}$ & 0.30 & $1.092\mathrm{E}{-1} \pm 7.811\mathrm{E}{-2}$ \\
            3 & 10 & $6.346\mathrm{E}{-3} \pm 1.163\mathrm{E}{-3}$ & 0.05 & $4.702\mathrm{E}{-2} \pm 1.237\mathrm{E}{-2}$ & 0.25 & $3.684\mathrm{E}{-1} \pm 8.453\mathrm{E}{-2}$ \\
            3 & 20 & $1.707\mathrm{E}{-2} \pm 1.161\mathrm{E}{-3}$ & 0.10 & $3.187\mathrm{E}{-2} \pm 1.407\mathrm{E}{-3}$ & 0.20 & $6.031\mathrm{E}{-1} \pm 4.876\mathrm{E}{-2}$ \\
            \addlinespace[2pt]
            4 & 5  & $2.415\mathrm{E}{-2} \pm 6.962\mathrm{E}{-3}$ & 0.60 & $1.008\mathrm{E}{-1} \pm 6.108\mathrm{E}{-3}$ & 0.40 & $9.324\mathrm{E}{-1} \pm 4.637\mathrm{E}{-2}$ \\
            4 & 10 & $9.084\mathrm{E}{-3} \pm 2.329\mathrm{E}{-3}$ & 0.40 & $4.960\mathrm{E}{-2} \pm 4.003\mathrm{E}{-3}$ & 0.35 & $9.405\mathrm{E}{-1} \pm 2.536\mathrm{E}{-2}$ \\
            4 & 20 & $1.684\mathrm{E}{-2} \pm 3.431\mathrm{E}{-4}$ & 0.05 & $2.732\mathrm{E}{-2} \pm 1.692\mathrm{E}{-3}$ & 0.25 & \textsc{OOM} \\
            \bottomrule
        \end{tabular}
    \end{table}

    \begin{figure}[htbp]
        \centering
        \includegraphics[width=0.99\linewidth]{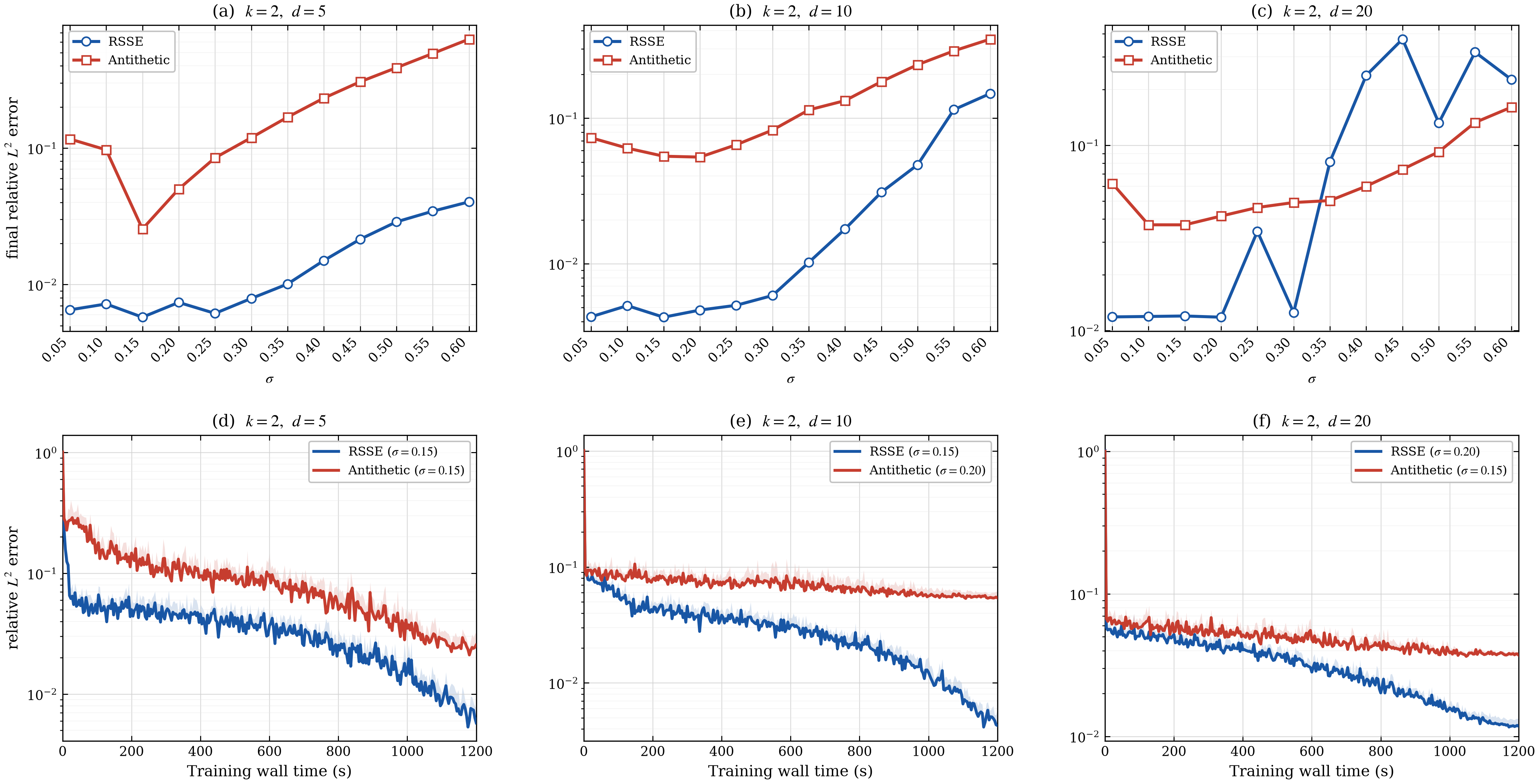}
        \caption{Numerical results for the polyharmonic PDE defined by \eqref{eq_poly} with $k = 2$.
            The top row reports the five-seed mean final $\operatorname{RE}$ defined in~\eqref{eq_rel_err} at each tested $\sigma$.
            The bottom row shows the test-time mean $\operatorname{RE}$ at the selected $\sigma$, with shading extending one sample standard deviation upward.}
        \label{fig_poly_k2}
    \end{figure}

    \begin{figure}[htbp]
        \centering
        \includegraphics[width=0.99\linewidth]{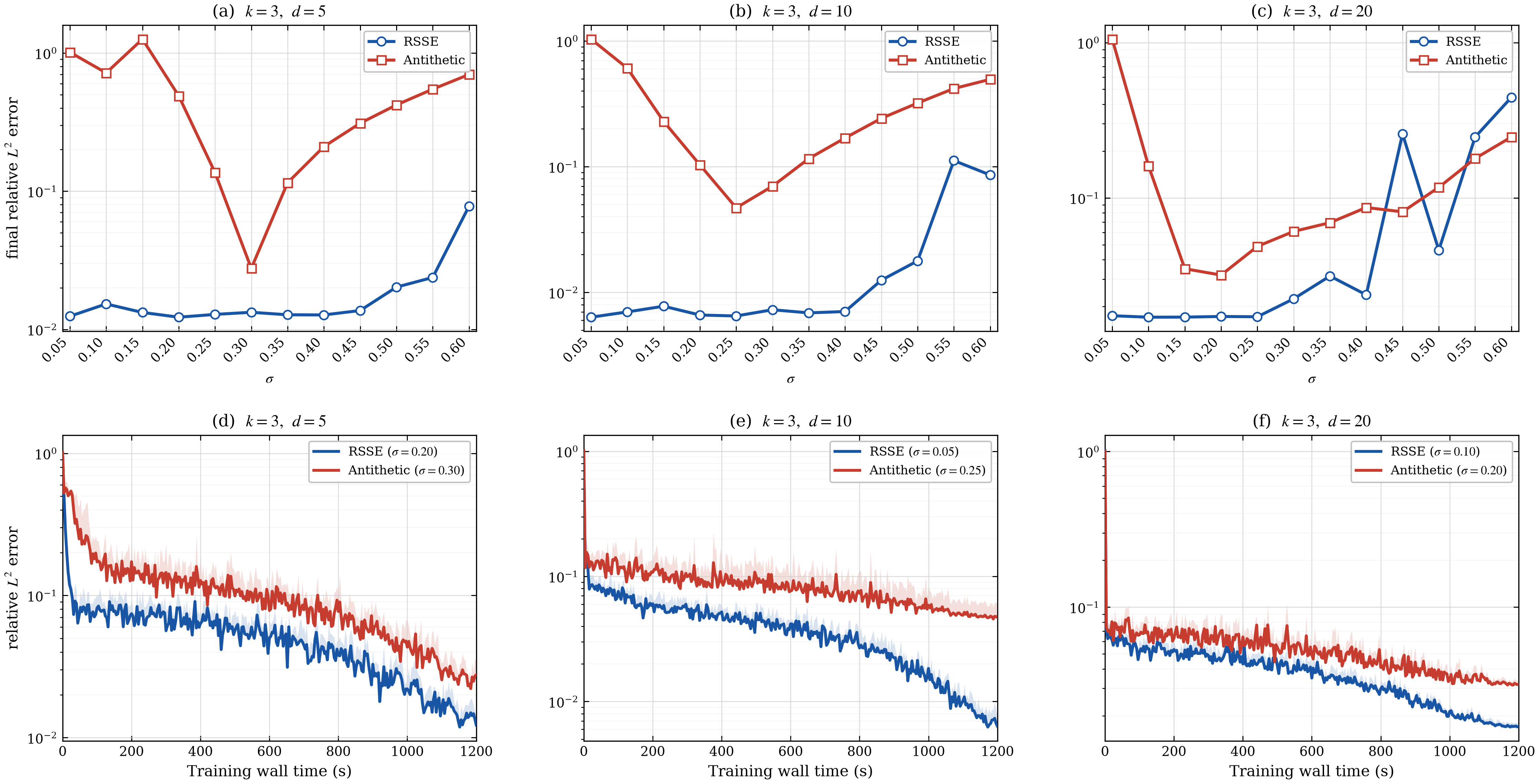}
        \caption{Numerical results for the polyharmonic PDE defined by \eqref{eq_poly} with $k = 3$.
            The top row reports the five-seed mean final $\operatorname{RE}$ defined in~\eqref{eq_rel_err} at each tested $\sigma$.
            The bottom row shows the test-time mean $\operatorname{RE}$ at the selected $\sigma$, with shading extending one sample standard deviation upward.}
        \label{fig_poly_k3}
    \end{figure}

    \begin{figure}[htbp]
        \centering
        \includegraphics[width=0.99\linewidth]{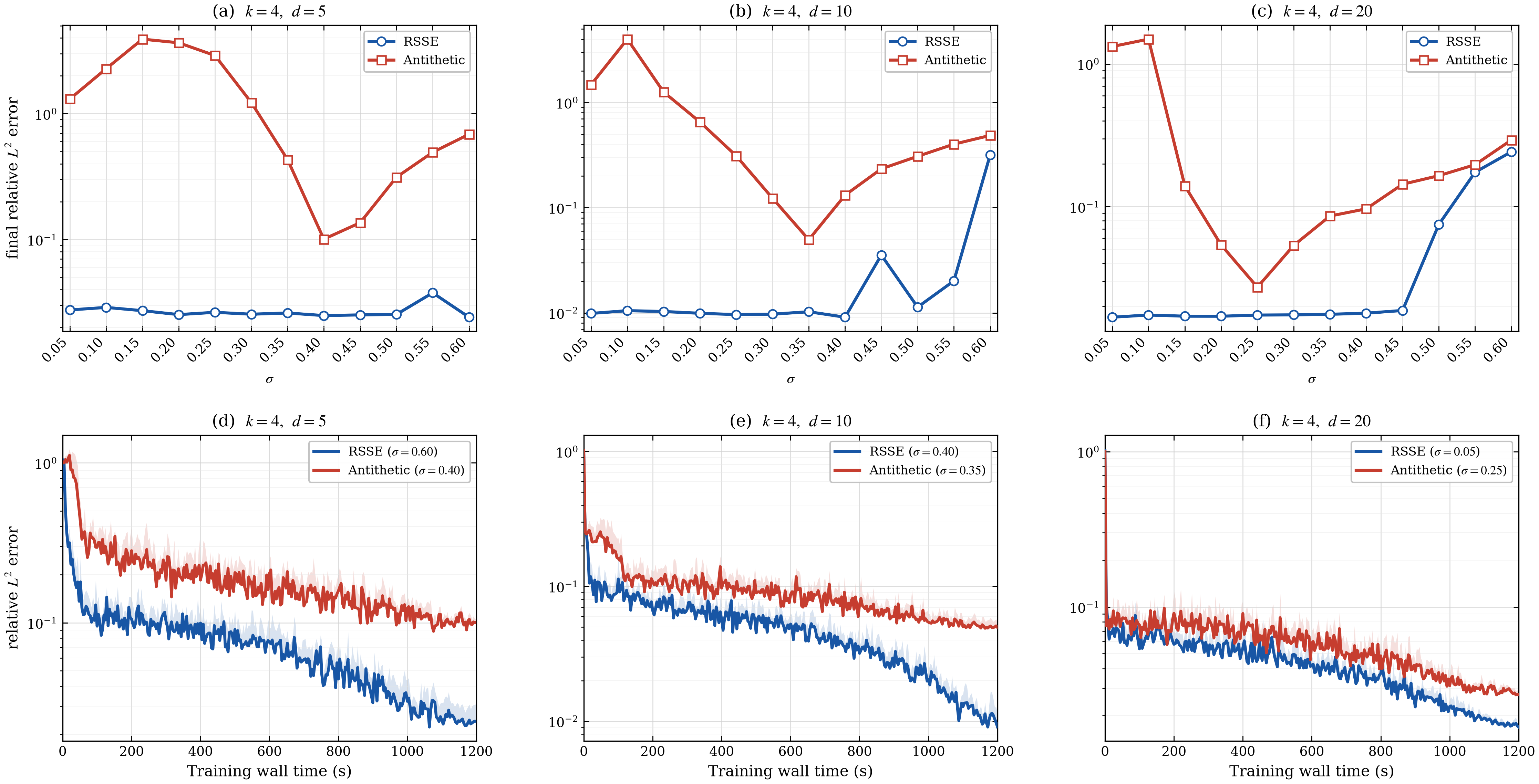}
        \caption{Numerical results for the polyharmonic PDE defined by \eqref{eq_poly} with $k = 4$.
            The top row reports the five-seed mean final $\operatorname{RE}$ defined in~\eqref{eq_rel_err} at each tested $\sigma$.
            The bottom row shows the test-time mean $\operatorname{RE}$ at the selected $\sigma$, with shading extending one sample standard deviation upward.}
        \label{fig_poly_k4}
    \end{figure}
    \FloatBarrier
\endgroup

\subsubsection{Fourth-order and sixth-order Cahn--Hilliard equations}\label{sec_ch}

For $m \in \{2, 3\}$ and $d \in \{5, 10, 20\}$, let $\Omega := [0, 1]^d$ and $\Gamma := \partial\Omega$.
With $\gamma_2 = 1$, $\gamma_{1, 2} = 10^{-2}$, and $\gamma_{1, 3} = -10^{-2}$, we solve the Cahn--Hilliard equation
\begin{equation}\label{eq_ch}
    \begin{cases}
        \partial_t u - \gamma_2\Delta(u^3 - u) + \gamma_{1, m}\Delta^m u = f_m(\bx, t), & (\bx, t) \in \Omega \times [0, 1], \\
        \partial_{\nu}^j u(\bx, t) = \partial_{\nu}^j u_*(\bx, t), & (\bx, t) \in \Gamma \times [0, 1], \quad j = 0, \ldots, m - 1, \\
        u(\bx, 0) = u_*(\bx, 0), & \bx \in \Omega,
    \end{cases}
\end{equation}
\par\noindent where $\nu$ denotes the outward unit normal.
The exact solution is
\begin{equation}\label{eq_ch_exact}
    u_*(\bx, t) := e^{-t}\exp\!\left(-\frac{|\bx|^2}{d}\right).
\end{equation}
The source term $f_m$ is obtained by substituting \eqref{eq_ch_exact} into the first line of \eqref{eq_ch}.
Unlike the hard constraints in Section~\ref{sec_poly}, the boundary and initial conditions here are imposed through squared mismatch terms in the training loss.

Table~\ref{tab_ch} reports the selected errors and step sizes, and Figures~\ref{fig_ch4}--\ref{fig_ch6} show the complete $\sigma$ sweeps and the corresponding training trajectories.
RSSE attains the smallest mean final $\operatorname{RE}$ in all six $(m, d)$ pairs.
The RSSE sweep curves are flatter over the tested step-size range, indicating lower sensitivity to $\sigma$.
AutoDiff is less accurate than RSSE throughout.

\begingroup
    \addtolength{\intextsep}{-4pt}
    \addtolength{\floatsep}{-4pt}

    \begin{table}[htbp]
        \centering
        \scriptsize
        \setlength{\tabcolsep}{2pt}
        \caption{Numerical results for the Cahn--Hilliard equation defined by \eqref{eq_ch}.
            We report the five-seed mean $ \pm $ sample standard deviation of the final $\operatorname{RE}$ defined in~\eqref{eq_rel_err}.
            For RSSE and antithetic sampling, $\sigma$ is selected to minimize the mean final $\operatorname{RE}$; AutoDiff is independent of $\sigma$.}
        \label{tab_ch}
        \begin{tabular}{ccccccc}
            \toprule
            $2m$ & $d$
            & \multicolumn{2}{c}{RSSE}
            & \multicolumn{2}{c}{antithetic}
            & \multicolumn{1}{c}{AutoDiff} \\
            \cmidrule(lr){3-4}\cmidrule(lr){5-6}\cmidrule(lr){7-7}
            & & \multicolumn{1}{c}{final $\operatorname{RE}$}
            & \multicolumn{1}{c}{$\sigma$}
            & \multicolumn{1}{c}{final $\operatorname{RE}$}
            & \multicolumn{1}{c}{$\sigma$}
            & \multicolumn{1}{c}{final $\operatorname{RE}$} \\
            \midrule
            4 & 5  & $1.567{\times}10^{-2} \pm 4.964{\times}10^{-3}$ & 0.30 & $3.114{\times}10^{-2} \pm 8.792{\times}10^{-3}$ & 0.20 & $3.862{\times}10^{-2} \pm 9.484{\times}10^{-3}$ \\
            4 & 10 & $1.632{\times}10^{-2} \pm 4.788{\times}10^{-3}$ & 0.20 & $2.628{\times}10^{-2} \pm 6.986{\times}10^{-3}$ & 0.30 & $4.951{\times}10^{-2} \pm 1.195{\times}10^{-2}$ \\
            4 & 20 & $1.991{\times}10^{-2} \pm 5.800{\times}10^{-3}$ & 0.40 & $3.565{\times}10^{-2} \pm 7.831{\times}10^{-3}$ & 0.30 & $6.084{\times}10^{-2} \pm 5.307{\times}10^{-3}$ \\
            \addlinespace[2pt]
            6 & 5  & $2.244{\times}10^{-2} \pm 1.252{\times}10^{-2}$ & 0.40 & $3.560{\times}10^{-2} \pm 3.427{\times}10^{-3}$ & 0.55 & $7.675{\times}10^{-2} \pm 4.383{\times}10^{-3}$ \\
            6 & 10 & $2.544{\times}10^{-2} \pm 9.976{\times}10^{-3}$ & 0.30 & $5.842{\times}10^{-2} \pm 1.437{\times}10^{-2}$ & 0.35 & $7.803{\times}10^{-2} \pm 3.071{\times}10^{-3}$ \\
            6 & 20 & $3.647{\times}10^{-2} \pm 1.717{\times}10^{-2}$ & 0.60 & $8.062{\times}10^{-2} \pm 1.678{\times}10^{-2}$ & 0.55 & $7.398{\times}10^{-2} \pm 5.364{\times}10^{-3}$ \\
            \bottomrule
        \end{tabular}
    \end{table}
    \FloatBarrier

    \begin{figure}[htbp]
        \centering
        \includegraphics[width=0.99\linewidth]{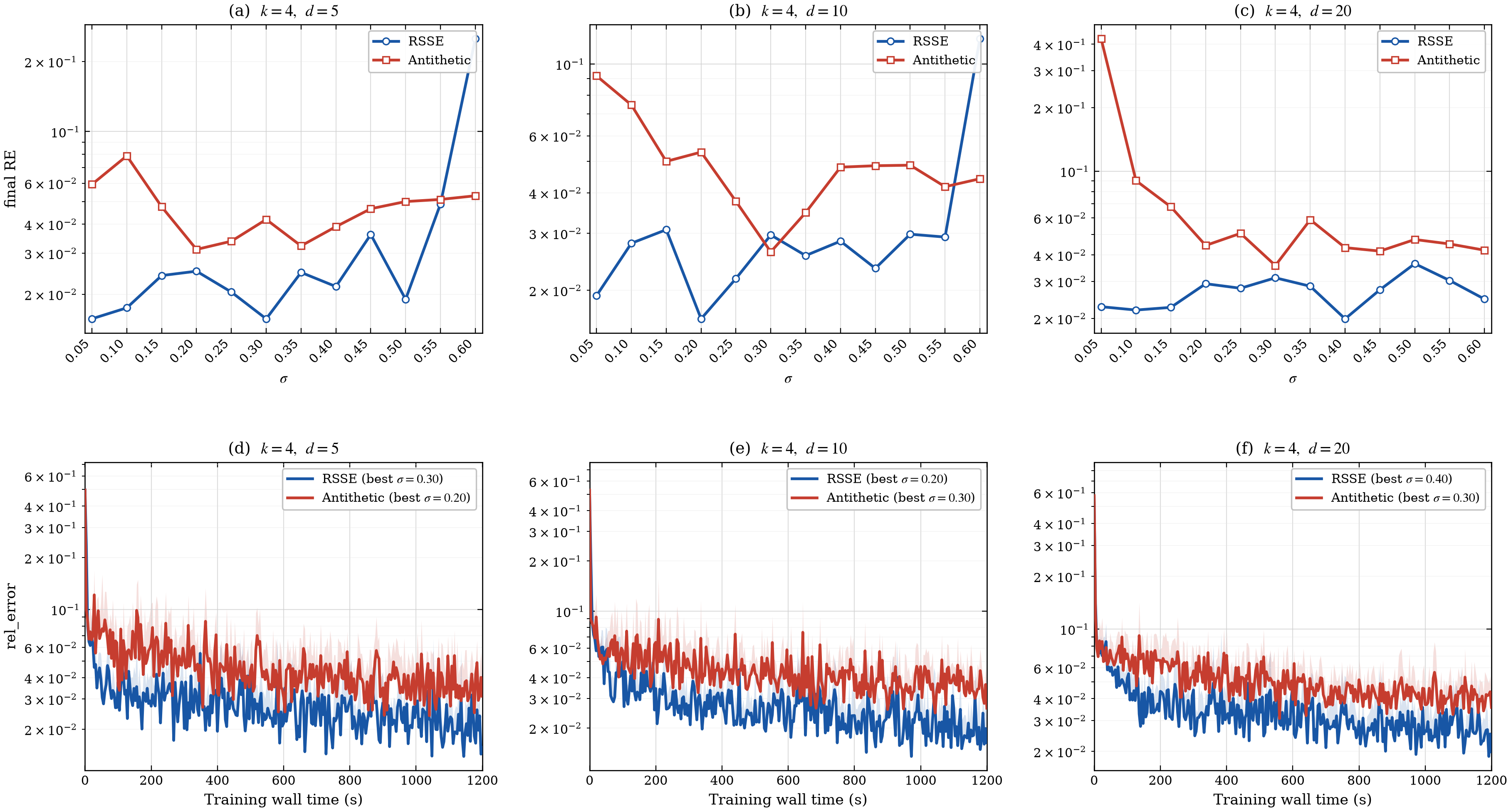}
        \caption{Numerical results for the Cahn--Hilliard equation defined by \eqref{eq_ch} with $m = 2$.
            In the panel titles, $k = 2m = 4$ denotes the differential order.
            The top row reports the five-seed mean final $\operatorname{RE}$ defined in~\eqref{eq_rel_err} at each tested $\sigma$.
            The bottom row shows the test-time mean $\operatorname{RE}$ at the selected $\sigma$, with shading extending one sample standard deviation upward.}
        \label{fig_ch4}
    \end{figure}

    \begin{figure}[htbp]
        \centering
        \includegraphics[width=0.99\linewidth]{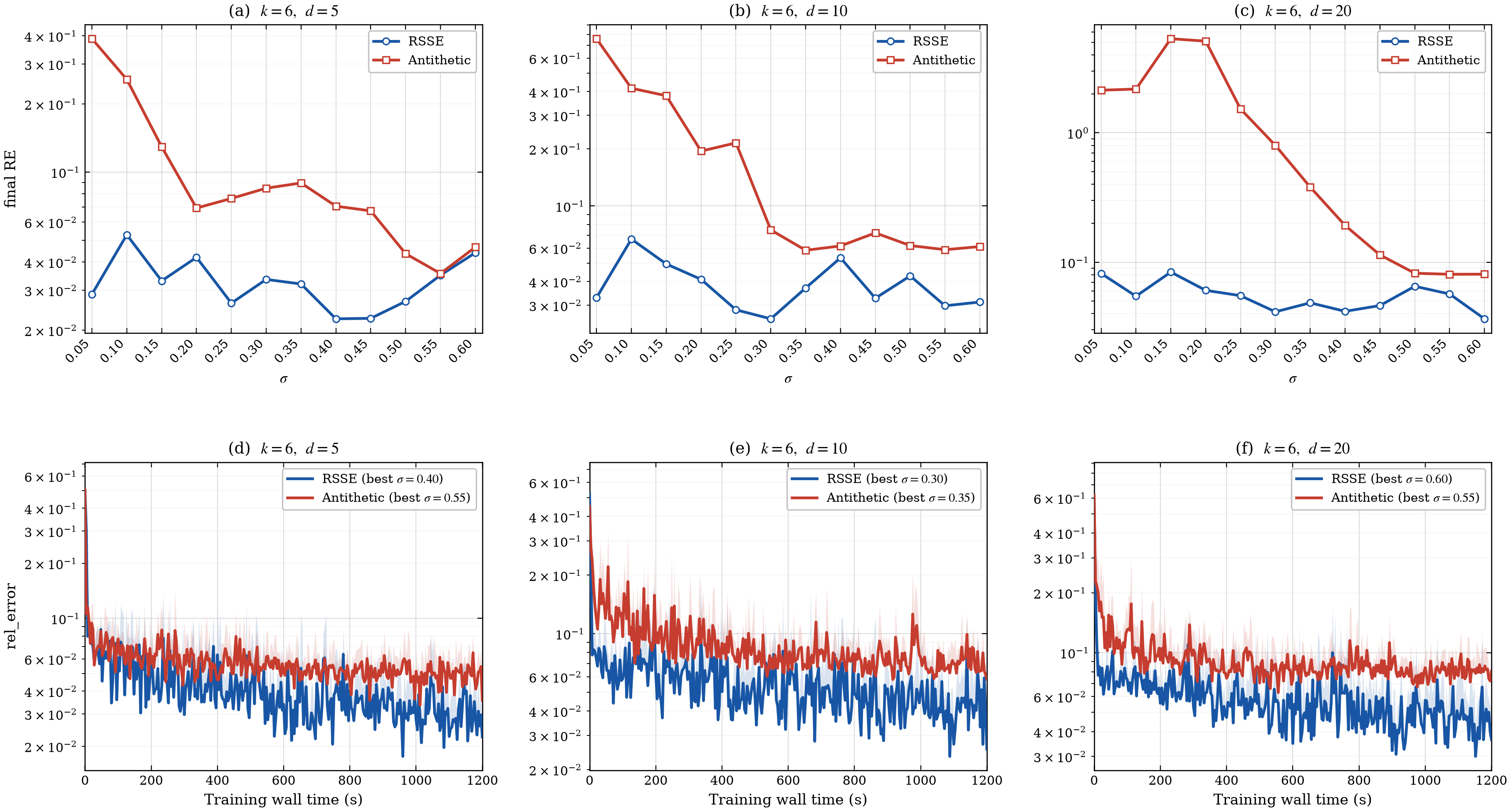}
        \caption{Numerical results for the Cahn--Hilliard equation defined by \eqref{eq_ch} with $m = 3$.
            In the panel titles, $k = 2m = 6$ denotes the differential order.
            The top row reports the five-seed mean final $\operatorname{RE}$ defined in~\eqref{eq_rel_err} at each tested $\sigma$.
            The bottom row shows the test-time mean $\operatorname{RE}$ at the selected $\sigma$, with shading extending one sample standard deviation upward.}
        \label{fig_ch6}
    \end{figure}
    \FloatBarrier
\endgroup


\section{Conclusion}\label{sec_concl}

We present a variance-stable RSSE for Monte Carlo approximation of high-order derivatives from function evaluations.
The method extends antithetic sampling from a symmetric pair on the real line to rotationally symmetric points in the complex plane, eliminating the lower-order contributions responsible for small-step variance divergence.
Under the stated assumptions, the estimator's variance converges to a finite limit as the step size vanishes, while its approximation bias is of order $\mathcal{O}(\sigma^n)$, where $\sigma$ is the step size and $n$ is the number of sampling points, chosen no smaller than the target derivative order.
We establish explicit bounds on the approximation bias and the deviation of the variance from its limiting value.
Numerical experiments support the predicted mean and variance expansions for second-, third-, and fourth-order derivative components.
The PDE experiments further illustrate the effectiveness of RSSE in PINNs for polyharmonic and Cahn--Hilliard equations under a fixed training-time budget.

\bibliographystyle{siam}
\bibliography{references}

\end{document}